# Online adaptive basis construction for nonlinear model reduction through local error optimization

Jun-Geol Ahn[1], Hyun-Ik Yang[1*], Jin-Gyun Kim[2*]

[1] *Department of Mechanical Design Engineering, Hanyang University, Wangsibri-ro 222, Seongdong-gu, Seoul, Republic of Korea*

[2] *Department of Mechanical Engineering (Integrated Engineering), Kyung Hee University, 1732 Deogyeong-daero, Giheung-gu, Yongin-si 17104, Republic of Korea*

**Abstract**

The accuracy of the reduced-order model (ROM) mainly depends on the selected basis. Therefore, it is essential to compute an appropriate basis with an efficient numerical procedure when applying ROM to nonlinear problems. In this paper, we propose an online adaptive basis technique to increase the quality of ROM while decreasing the computational costs in nonlinear problems. In the proposed method, the adaptive basis is defined by the low-rank update formulation, and two auxiliary vectors are set to implement this low-rank condition. To simultaneously tackle the issues of accuracy and the computational cost of the ROM basis, the auxiliary vectors are algebraically derived by optimizing a local residual operator. As a result, the reliability of ROM is significantly improved with a low computational cost because the error information can be contained without inverse operations of the full model dimension required in conventional approaches. The other feature of the proposed iterative algorithm is that the number of the initial incremental ROM basis could be varied, unlike in the typical online adaptive basis approaches. It may provide a fast and effective spanning process of the high-quality ROM subspace in the iteration step. A detailed derivation process of the proposed method is presented, and its performance is evaluated in various nonlinear numerical examples.

**Keywords**: model order reduction; online adaptive basis; nonlinear dynamic analysis; adaptive strategy; proper orthogonal decomposition

* Corresponding authors: Jin-Gyun Kim (jingyun.kim@khu.ac.kr) and Hyun-Ik Yang (skynet@hanyang.ac.kr).

# 1. Introduction

The projection-based reduced-order model (ROM) has become an essential tool to effectively solve a large and complex full-order model (FOM) [1-5]. In linear problems, various ROM techniques have been widely used to obtain a reliable basis and results. The excellent performances of the developed ROM have been verified in many engineering practices [6-8]. However, those ROM techniques are limited in representing nonlinear systems.

A proper orthogonal decomposition (POD) is a popular ROM method for nonlinear problems [9-12]. The orthonormal basis returned by the POD method can effectively compress interesting data in certain problems. Therefore, the POD method has been employed not only in the ROM approaches but also in many computational practices, such as neural networks [13-14].

However, if data corrected by the POD method are insufficient to represent important features of the problems, POD-based ROM may lead to untrustworthy results. To solve this problem, online adaptive basis techniques have been developed [15,17-18]. In these methods, the additional basis is newly obtained when ROM cannot satisfy an error tolerance during the online calculation stage. As a result, a modified subspace with the additional ROM basis becomes more appropriate for representing the nonlinear response at a certain step than the initial subspace. Thus, a reconstructed ROM can provide better solution accuracy [17-19]. Despite an excellent performance of existing online adaptive basis techniques, several challenges remain. For instance, adding an additional basis can improve the ROM quality because of an extended subspace of ROM [15-16]. However, a huge number of basis vectors may be required because the convergence rate from the ROM to FOM solutions depends not only on the number of basis vectors. A residual operator that represents the error information is used to solve the convergence rate problem [17-18]. This is because the additional basis obtained from the error term can directly reduce the differences between the ROM and FOM responses. Nevertheless, it may incur relatively high computational costs during the online process because some numerical procedures for the residual operator include inverse operations in the FOM dimension [15,20].

In this paper, we propose a new online adaptive basis technique to ensure computational efficiency when considering the error information. In the proposed method, the additional basis is defined from the low-rank update formulation and its auxiliary vectors. These auxiliary vectors are algebraically derived from the optimal solutions of a local residual operator to reduce the calculation costs while including the error term. Moreover, to systematically guarantee the existence of the derived auxiliary vectors, we employ and transform the proven lemma that gives the conditions for reducing the Frobenius norm of the matrix (or vector) [19-20]. Importantly, the local residual operator can be constructed by freely selecting the components of the residual operator. Thus, the inverse operations of the full model dimension required in conventional approaches can be avoided in the suggested approach [15,17-18]. Consequently, the proposed method can have robust existence conditions and can significantly enhance the ROM accuracy with a low calculation cost. Furthermore, the provided algorithm can also be employed to quickly expand the high-quality ROM subspace. This is because the additional basis in the proposed method is implemented with a group of vectors using the low-rank update condition, unlike in the existing methods that return the additional basis as a vector.

The remainder of this paper is organized as follows. Section 2 reviews the problem setup for a nonlinear problem and the online adaptive basis technique. The proposed method is derived and explained in Section 3. Section 4 shows the performance of the proposed method in various nonlinear problems. Conclusions are given in Section 5.

## 2. Problem setup for the nonlinear problem and the online adaptive basis technique.

In this section, we formulate a nonlinear problem and its projection-based ROM. In addition, we review the POD method and the online adaptive basis procedure [9-11, 15-18].

*2.1. Full-order model and its projection-based ROM*

A parametrized nonlinear dynamic problem can be represented as follows [15-16]:

$$\dot{x} = f(x; \mu), \tag{1a}$$
$$x \in \mathbb{R}^N, \ f \in \mathbb{R}^N, \ \mu \subset D \in \mathbb{R}^{N_\mu}, \tag{1b}$$

where $x$ indicates the solution, and $f$ denotes the nonlinear expression. Input parameters and their parameter domain are expressed as $\mu$ and $D$, respectively. Superscripts $N$ and $N_\mu$ represent the dimensions of the solution and input parameters, respectively. In Eq. (1), the first-order differentiation of time for solution $x$ is defined as $\dot{x}$ (i.e., $\dot{x} := \frac{d}{dt}x$).

The nonlinear problem in Eq. (1) or its solver formulation can be written using a residual operator as follows:

$$\boldsymbol{r}_k(\boldsymbol{x}_k; \boldsymbol{\mu}) = \boldsymbol{0}, \quad \boldsymbol{r}_k: \mathbb{R}^N \times \mathbb{R}^{N_u} \to \mathbb{R}^N, \tag{2}$$

where $\boldsymbol{r}_k$ denotes the residual operator at iteration $k$. Notably, various states of the nonlinear problem or a discretized solver formulation can be represented using Eq. (2) [15-16]. For example, the nonlinear problem in Eq. (1) can be formulated using the residual operator $\boldsymbol{r}_k$ as follows:

$$\boldsymbol{r}_k = \dot{\boldsymbol{x}}_k - \boldsymbol{f}(\boldsymbol{x}_k; \boldsymbol{\mu}). \tag{3}$$

Moreover, the implicit backward Euler method for Eq. (1) can be described in the residual operator form:

$$\boldsymbol{r}_k = (\boldsymbol{x}_k - \boldsymbol{x}_{k-1}) - \Delta t \cdot \boldsymbol{f}(\boldsymbol{x}_k; \boldsymbol{\mu}), \tag{4}$$

where $\Delta t$ denotes the discretized time step. The residual operator $\boldsymbol{r}_k$ in Eq. (2) is a more general form of that in Eq. (1), see Eqs. (1) to (4). Therefore, in this work, we use Eq. (2) instead of Eq. (1). To improve readability, Eq. (2) is called a full model.

The projection-based ROM for the full model in Eq. (2) is written as follows [4-8]:

$$\boldsymbol{W}^T \boldsymbol{r}_k(\hat{\boldsymbol{x}}_k; \boldsymbol{\mu}) = \boldsymbol{0}, \tag{5a}$$
$$\boldsymbol{x}_k \approx \hat{\boldsymbol{x}}_k = \bar{\boldsymbol{x}} + \boldsymbol{\Phi} \boldsymbol{q}_k, \tag{5b}$$
$$\boldsymbol{W} \in \mathbb{R}^{N \times n}, \; \boldsymbol{\Phi} \in \mathbb{R}^{N \times n}, \; \hat{\boldsymbol{x}} \in \mathbb{R}^N, \; \bar{\boldsymbol{x}} \in \mathbb{R}^N, \; \boldsymbol{q}_k \in \mathbb{R}^n, \tag{5c}$$

where $\boldsymbol{W}$ and $\boldsymbol{\Phi}$ denote the left subspace and the trial basis for ROM, respectively. Moreover, $\hat{\boldsymbol{x}}$ and $\boldsymbol{q}$ in Eq. (5) indicate an approximated solution and the ROM solution, respectively. In Eq. (5), the reference configuration (e.g., the initial guess) is expressed by $\bar{x}$. Superscript $n$ indicates the selected basis size $n$, which can be much smaller than the full model size $N$ (i.e., $n \ll N$). Thus, computational costs can be significantly reduced by using Eq. (5) instead of the full model in Eq. (2). Note that Eq. (5) is called Galerkin projection or Petrov–Galerkin projection, depending on the left subspace $W$ [21]. For instance, Galerkin projection for Eq. (5) is determined as

$$\boldsymbol{\Phi}^T \boldsymbol{r}_k(\bar{x} + \boldsymbol{\Phi} \boldsymbol{q}_k; \boldsymbol{\mu}) = \boldsymbol{0}. \tag{6}$$

In the projection-based ROM, the trial basis $\boldsymbol{\Phi}$ in Eqs. (5) and (6) can be freely taken by users. Nevertheless, the quality of ROM highly depends on the trial basis $\boldsymbol{\Phi}$. Therefore, a selection of the trial basis $\boldsymbol{\Phi}$ is an important issue in ROM. In linear problems, the trial basis $\boldsymbol{\Phi}$ is usually chosen using the generalized eigenvalue problem of system matrices at the offline process [6-7]. This is because the linear model and its dynamic characteristics are invariant during the online calculation. However, because the nonlinear model and its properties may change in the online stage, the generalized eigenpairs obtained from the offline stage may not effectively explain the responses. For this reason, the POD method is employed to acquire the trial basis $\boldsymbol{\Phi}$ in nonlinear problems [9-12]. In the next section, we briefly introduce the POD method.

*2.2. Proper orthogonal decomposition (POD)*

POD is a popular method to obtain the ROM basis in nonlinear problems. To compute the basis, the POD method requires a collection of solution vectors, called the snapshot matrix. The snapshot matrix can be constructed as follows [9-12]:

$$\boldsymbol{X} = [\boldsymbol{x}(\boldsymbol{\mu}_1) \quad \cdots \quad \boldsymbol{x}(\boldsymbol{\mu}_M)], \ \boldsymbol{X} \in \mathbb{R}^{N \times M}, \tag{7}$$

where $\boldsymbol{X}$ denotes the snapshot matrix, and superscript $M$ indicates the number of solutions in Eq. (7). The basis of POD method is computed from the snapshot matrix $\boldsymbol{X}$. Thus, a selection of solutions has a significant impact to the POD-based ROM.

The POD method finds an orthonormal basis using the following minimization problem [9-12]:

$$\min_{v_1,\cdots,v_p} \sum_{i=1}^{M}\left(\left\|x(\mu_i)-\sum_{j=1}^{p}\left(v_j^T x(\mu_i)\right)v_j\right\|_2\right)^2, \tag{8a}$$

$$V=[v_1 \quad \cdots \quad v_p]\in\mathbb{R}^{N\times p},\ V\in\mathbb{R}^{N\times p}, \tag{8b}$$

where $V$ denotes the orthonormal basis acquired from the minimization problem in Eq. (8). Moreover, $\|\cdot\|_2$ and superscript $T$ represent a second norm and transpose operators of the matrix (or vector), respectively. The number of basis vectors in $V$ is shown as superscript $p$.

As a result, the trial basis $\Phi$ of POD-based ROM is implemented by using basis $V$ in Eq. (8) as follows:

$$\Phi=[v_1 \quad \cdots \quad v_n]\in\mathbb{R}^{N\times n},\ n\leq p. \tag{9}$$

The trial basis $\Phi$ in Eq. (9) has shown excellent performance in various nonlinear problems [9-14]. However, it is not easy to confirm that the constructed snapshot matrix $X$ sufficiently represents the required nonlinear problems. Importantly, if the snapshot matrix does not contain enough information, the reliability of POD-based ROM may significantly decrease in the online stage. To address the quality issue of the ROM basis, the online adaptive basis procedure has been developed [15,17-18], which is briefly discussed in the next section.

*2.3. Online adaptive basis technique for solving nonlinear problems*

In the online adaptive basis procedure, the trial basis $\Phi$ is recomputed when ROM does not converge during the online calculation. In particular, adding the other basis to the trial basis $\Phi$ can easily increase the accuracy of ROM. An extended result of the trial basis $\Phi$ can be written as follows:

$$\Phi \leftarrow [\Phi \quad \Psi],\ \Psi\in\mathbb{R}^{N\times u}, \tag{10}$$

where $\Psi$ denotes an additional basis, and the size of $\Psi$ is represented by superscript $u$. A subspace of the updated basis $[\Phi \quad \Psi]$ is expanded compared with the subspace of the trial basis $\Phi$ [17-18]. As a result, the updated basis $[\Phi \quad \Psi]$ can easily increase the fidelity of ROM [15-16]. Thus, in this study, we focus on determining the additional basis $\boldsymbol{\Psi}$. Algorithm 1 describes the online adaptive basis technique using Galerkin projection in Eq. (6) [2, 14-16, 21]. Note that Algorithm 1 can be also applied to Petrov–Galerkin projection with the appropriate left subspaces $\boldsymbol{W}$ [21].

---

Algorithm 1. Online adaptive basis technique using Galerkin projection

---

Input: iteration $\boldsymbol{k}$, trial basis $\boldsymbol{\Phi}$, ROM solver tolerance $\varepsilon_{\text{ROM}}$, and FOM solver tolerance $\varepsilon_{\text{FOM}}$.

Output: 1) Updated basis $[\boldsymbol{\Phi} \quad \boldsymbol{\Psi}]$ and 2) ROM solution $\boldsymbol{q}_k$.

1. Compute ROM solution $\boldsymbol{q}_k$ satisfying ROM solver tolerance $\varepsilon_{\text{ROM}}$:

    $\|\boldsymbol{\Phi}^T \boldsymbol{r}_k(\bar{\boldsymbol{x}} + \boldsymbol{\Phi}\boldsymbol{q}_k; \boldsymbol{\mu})\|_2 \leq \varepsilon_{\text{ROM}}$.

2. Calculate FOM error $\varepsilon$:

    $\varepsilon = \|\boldsymbol{r}_k(\bar{\boldsymbol{x}} + \boldsymbol{\Phi}\boldsymbol{q}_k; \boldsymbol{\mu})\|_2$.

3. If error $\varepsilon$ in Step 2 does not converge (i.e., $\varepsilon > \varepsilon_{\text{FOM}}$) do

4.     Take the additional basis $\boldsymbol{\Psi}$.

5.     Construct the updated basis $[\boldsymbol{\Phi} \quad \boldsymbol{\Psi}]$.

6.     Replace the trial basis $\boldsymbol{\Phi}$ using the updated basis $[\boldsymbol{\Phi} \quad \boldsymbol{\Psi}]$.

    $\boldsymbol{\Phi} \leftarrow [\boldsymbol{\Phi} \quad \boldsymbol{\Psi}]$.

7.     Return to Step 1.

8. end if.

---

There are several ways to take the additional basis $\boldsymbol{\Psi}$. For example, if the trial basis $\boldsymbol{\Phi}$ is obtained using the POD basis in Eq. (9), the additional basis $\boldsymbol{\Psi}$ can be easily chosen as the remaining POD basis in Eq. (8) as follows:

$$\boldsymbol{\Psi} = [\boldsymbol{v}_{n+1} \quad \cdots \quad \boldsymbol{v}_{n+u}], \quad n + u \leq p. \tag{11}$$

In this study, when the additional basis $\boldsymbol{\Psi}$ is determined using Eq. (11), it is called POD-ROM for convenience.

Researches of various nonlinear problems have confirmed that using the additional basis $\boldsymbol{\Psi}$ increases the ROM accuracy [15-18]. Nevertheless, the quality of the ROM basis does not entirely depend on the number of basis vectors. Therefore, some improvements are required. For instance, if the updated basis in Algorithm 1 are derived without the error information, the differences between the ROM and FOM responses may not sufficiently decrease, even if the ROM subspace is expanded. To directly tackle the problem of differences during the online process, an effect of the error $\varepsilon$ and its optimal condition should be considered in the adaptive procedure [17-18]. In this study, we cover the Newton method and Galerkin projection to explicitly demonstrate the adaptive process that involves error information.

Using the Newton method, the residual operator $r_k$ in Eq. (2) can be written as follows:

$$r_k = J_k \Phi q_k + F_k, \tag{12a}$$

$$F_k = F_k(\bar{x} + \Phi q_k; \mu), \quad J_k = J_k(\bar{x} + \Phi q_k; \mu), \tag{12b}$$

$$F_k \in \mathbb{R}^N, \quad J_k \in \mathbb{R}^{N \times N}, \tag{12c}$$

where $F_k$ and $J_k$ denote the nonlinear function and the Jacobian matrix in the Newton method, respectively. Algorithm 2 describes the online adaptive basis technique that considers the effect of the error $\varepsilon$ using the Newton method and Galerkin projection [17-18, 21].

| Algorithm 2. Online adaptive basis technique for the Newton method and Galerkin projection that considers the error information. |
|---|

Input: iteration $k$, trial basis $\boldsymbol{\Phi}$, ROM solver tolerance $\varepsilon_{\text{ROM}}$, and FOM solver tolerance $\varepsilon_{\text{FOM}}$.

Output: 1) Updated basis $[\boldsymbol{\Phi} \quad \boldsymbol{\Psi}]$ and 2) ROM solution $\boldsymbol{q}_k$.

1. Compute the ROM solution $\boldsymbol{q}_k$ using the Newton method:

   $\|\boldsymbol{\Phi}^T \boldsymbol{r}_k\|_2 \leq \varepsilon_{\text{ROM}}, \ \boldsymbol{r}_k = \boldsymbol{J}_k \boldsymbol{\Phi} \boldsymbol{q}_k + \boldsymbol{F}_k$.

2. Calculate the FOM error $\varepsilon$:

   $\varepsilon = \|\boldsymbol{r}_k\|_2$.

3. If error $\varepsilon$ in Step 2 does not converge (i.e., $\varepsilon > \varepsilon_{\text{FOM}}$) do

4.     Choose the additional basis $\boldsymbol{\Psi}$ to satisfy the following minimization problem:

   $\boldsymbol{\Psi} = \underset{\widetilde{\boldsymbol{\Psi}}}{argmin} \left( \|\boldsymbol{J}_k \widetilde{\boldsymbol{\Psi}} - \boldsymbol{r}_k\|_F \right)^2$.

5.     Construct the updated basis $[\boldsymbol{\Phi} \quad \boldsymbol{\Psi}]$.

6.     Replace basis $\boldsymbol{\Phi}$ using the updated basis $[\boldsymbol{\Phi} \quad \boldsymbol{\Psi}]$:

   $\boldsymbol{\Phi} \leftarrow [\boldsymbol{\Phi} \quad \boldsymbol{\Psi}]$.

7.     Return to Step 1.

8. end if

In Algorithm 2, $\|\cdot\|_F$ denotes the Frobenius norm operator of the matrix (or vector). Note that the thin-QR factorization can be applied in Algorithm 2 to make an orthonormal condition [20]. In addition, the initialization of the basis can be also used to limit the size of the updated basis $[\boldsymbol{\Phi} \quad \boldsymbol{\Psi}]$ [15-16].

The additional basis $\boldsymbol{\Psi}$ calculated from Algorithm 2 can consider error $\varepsilon$. Thus, differences between the ROM and FOM responses can be immediately decreased [17-18]. However, in this case, Algorithm 2 requires the inverse operation at the FOM dimension. Therefore, the online process may incur huge computational costs [15, 20]. In this study, when the additional basis $\boldsymbol{\Psi}$ is obtained using Algorithm 2, it is called F-ROM for convenience.

In order to take the error effect with efficient computational costs, we suggest a new adaptive procedure using the low-rank update formulation and an optimal condition of a local residual operator. The proposed method is derived and illustrated in the following section.

# 3. Proposed method for the online adaptive basis technique

The additional basis $\boldsymbol{\Psi}$ in the proposed method is derived in a three-step process. First, when the additional basis $\boldsymbol{\Psi}$ is required because of error $\varepsilon$ (Step 3 in Algorithm 1), a local residual operator is defined by selecting only a portion of the residual operator $\boldsymbol{r}_k$. Second, two auxiliary vectors are derived by optimizing the local residual operator. Third, the additional basis $\boldsymbol{\Psi}$ is calculated using the trial basis and the low-rank update formulation. Here, the low-rank condition is represented using auxiliary vectors at the second step. As a result, in the proposed method, the additional basis $\boldsymbol{\Psi}$ can effectively compensate for the error effects without the inverse operation in the FOM dimension, described in Section 2.

To explicitly describe the proposed method, we employ the Newton method and Galerkin projection in this study. Importantly, the proposed approach can also be applied to other solvers and/or projections [1-2, 6-9, 21].

The local residual operator in the Newton method is determined using the Boolean matrix and residual operator $\boldsymbol{r}_k$ as follows:

$$\boldsymbol{L}_r = \boldsymbol{B}_k^T \boldsymbol{r}_k = \boldsymbol{B}_k^T (\boldsymbol{J}_k \boldsymbol{\Phi} \boldsymbol{q}_k + \boldsymbol{F}_k), \tag{13a}$$

$$\boldsymbol{L}_r \in \mathbb{R}^{n_{sel}}, \boldsymbol{B}_k \in \mathbb{R}^{N \times n_{sel}}, \tag{13b}$$

where $\boldsymbol{L}_r$ and $\boldsymbol{B}_k$ denote the local residual operator and the Boolean matrix, respectively. The size of $\boldsymbol{L}_r$ is represented by superscript $n_{sel}$. Here, $\boldsymbol{B}_k$ is used to collect the components of the residual operator $\boldsymbol{r}_k$. For example, if users want to obtain the first and third rows using $\boldsymbol{B}_k \in \mathbb{R}^{3 \times 2}$, the Boolean matrix is constructed as $\boldsymbol{B}_k = [\boldsymbol{e}_1 \quad \boldsymbol{e}_3]$, where $\boldsymbol{e}_1 = [1 \quad 0 \quad 0]^T$ and $\boldsymbol{e}_1 = [0 \quad 0 \quad 1]^T$. As the elements of the residual operator $\boldsymbol{r}_k$ can be freely decided by users, the size of the local residual operator $\boldsymbol{L}_r$ can be significantly reduced compared with the full model size (i.e., $n_{sel} \ll N$).

Unlike the existing approaches in Algorithm 2, the additional basis $\boldsymbol{\Psi}$ is formulated using the trial basis $\boldsymbol{\Phi}$ and the low-rank update formulation. Thus, the proposed additional basis $\boldsymbol{\Psi}$ is represented as follows:

$$\boldsymbol{\Psi} = \boldsymbol{\Phi} + \boldsymbol{\xi}\boldsymbol{\psi}^T, \tag{14a}$$

$$\boldsymbol{\xi} \in \mathbb{R}^N, \; \boldsymbol{\psi} \in \mathbb{R}^n, \tag{14b}$$

where $\boldsymbol{\xi}$ and $\boldsymbol{\psi}$ denote the auxiliary vectors in the proposed method. Here, the rank of $\boldsymbol{\xi}\boldsymbol{\psi}^T$ is 1 because all vectors in $\boldsymbol{\xi}\boldsymbol{\psi}^T$ can be represented as $[\psi_1 \boldsymbol{\xi} \; \cdots \; \psi_n \boldsymbol{\xi}]$, where $\boldsymbol{\psi}^T = [\psi_1 \; \cdots \; \psi_n]$. For this reason, the $\boldsymbol{\xi}\boldsymbol{\psi}^T$ form is called as the low-rank update formulation.

In this work, the purpose of the additional basis $\boldsymbol{\Psi}$ is to reduce the quantity of the local residual operator $\boldsymbol{L}_r$. In addition, the quality of the additional basis $\boldsymbol{\Psi}$ can be improved by considering the Jacobian matrix $\boldsymbol{J}_k$ [17-18], as can be verified using conventional methods, see Step 4 in Algorithm 2. Therefore, we derive the additional basis $\boldsymbol{\Psi}$ that decreases the norm of the local residual operator $\boldsymbol{L}_r$ while reflecting the effects of the Jacobian matrix $\boldsymbol{J}_k$. To explain this objective, Eqs. (13) and (14) are combined. Substituting Eq. (14) into the location of the trial basis $\boldsymbol{\Phi}$ in Eq. (13), the modified formulation of Eq. (13) is obtained as follows:

$$\boldsymbol{B}_k^T(\boldsymbol{J}_k \boldsymbol{\Psi} \boldsymbol{q}_k + \boldsymbol{F}_k) = \boldsymbol{L}_r + \boldsymbol{L}_J \boldsymbol{\xi} \boldsymbol{\psi}^T \boldsymbol{q}_k, \tag{15a}$$

$$\boldsymbol{L}_J = \boldsymbol{B}_k^T \boldsymbol{J}_k, \; \boldsymbol{L}_J \in \mathbb{R}^{n_{sel} \times N}, \tag{15b}$$

where $\boldsymbol{L}_J$ indicates the local Jacobian matrix. Note that the additional basis $\boldsymbol{\Psi}$ in Eq. (15) can demonstrate the influence of Jacobian matrix $\boldsymbol{J}_k$ using the local Jacobian matrix $\boldsymbol{L}_J$.

Using Eq. (15), the aim of the proposed method can be explicitly written as follows:

$$\left(\left\|\boldsymbol{L}_r + \boldsymbol{L}_J \boldsymbol{\xi} \boldsymbol{\psi}^T \boldsymbol{q}_k\right\|_F\right)^2 < (\|\boldsymbol{L}_r\|_F)^2. \tag{16}$$

Auxiliary vectors $\boldsymbol{\xi}$ and $\boldsymbol{\psi}$ in Eq. (16) should be obtained to calculate Eq. (14). However, the existence of these auxiliary vectors is not identified yet. Thus, it may be impossible to calculate

Eq. (14). To guarantee the presence of $\boldsymbol{\xi}$ and $\boldsymbol{\psi}$ in Eq. (16), we employ a previously proven lemma:

Lemma 1. Assume that $\boldsymbol{A} \in \mathbb{R}^{a \times c}$ and $\boldsymbol{C} \in \mathbb{R}^{b \times c}$ are arbitrary real matrices, and the rank of $\boldsymbol{C}$ is the number of rows in $\boldsymbol{C}$ (i.e., rank($\boldsymbol{C}$) = b). Then, there exist real variables $\boldsymbol{\alpha} \in \mathbb{R}^a$ and $\boldsymbol{\beta} \in \mathbb{R}^b$ with $(\|\boldsymbol{A} + \boldsymbol{\alpha}\boldsymbol{\beta}^T \boldsymbol{C}\|_F)^2 < (\|\boldsymbol{A}\|_F)^2$ if and only if $\|\boldsymbol{A}\boldsymbol{C}^T\|_F > 0$.

*Proof.* See reference [19]. □

Lemma 1 represents the conditions when the Frobenius norm of the matrix (or vector) can be decreased. Nevertheless, Lemma 1 is not directly applicable to Eq. (16) because the rank of the ROM solution $\boldsymbol{q}_k$ is not equal to the number of rows in $\boldsymbol{q}_k$. Therefore, to employ Lemma 1, we firstly use the following equation:

$$\boldsymbol{L}_r + \hat{\boldsymbol{\xi}}\hat{\boldsymbol{\psi}}^T d_k, \quad d_k = \|\boldsymbol{q}_k\|_2, \tag{17}$$

where $\hat{\boldsymbol{\xi}}$ and $\hat{\boldsymbol{\psi}}$ are the modified auxiliary variables, and $d_k$ denotes the magnitude of the ROM solution $\boldsymbol{q}_k$. The minimization problem for Eq. (17) can be expressed as follows:

$$\tilde{\boldsymbol{\xi}}, \tilde{\boldsymbol{\psi}} = \underset{\hat{\boldsymbol{\xi}},\hat{\boldsymbol{\psi}}}{\text{argmin}} \left(\left\|\boldsymbol{L}_r + \hat{\boldsymbol{\xi}}\hat{\boldsymbol{\psi}}^T d_k\right\|_F\right)^2, \tag{18a}$$

$$\left(\left\|\boldsymbol{L}_r + \tilde{\boldsymbol{\xi}}\tilde{\boldsymbol{\psi}}^T d_k\right\|_F\right)^2 < (\|\boldsymbol{L}_r\|_F)^2, \tag{18b}$$

where $\tilde{\boldsymbol{\xi}}$ and $\tilde{\boldsymbol{\psi}}$ indicate the optimal solutions for the modified auxiliary variables $\hat{\boldsymbol{\xi}}$ and $\hat{\boldsymbol{\psi}}$, respectively. Here, we want to show that unlike Eq. (16), Lemma 1 can be applied to Eq. (18). Using Lemma 1, the following Lemma for Eq. (18) can be obtained:

Lemma 2. There exist two variables $\tilde{\boldsymbol{\xi}}$ and $\tilde{\boldsymbol{\psi}}$ with $\left(\left\|\boldsymbol{L}_r + \tilde{\boldsymbol{\xi}}\tilde{\boldsymbol{\psi}}^T d_k\right\|_F\right)^2 < (\|\boldsymbol{L}_r\|_F)^2$ unless all values of $\boldsymbol{L}_r$ and $\boldsymbol{q}_k$ are zero.

*Proof.* $\boldsymbol{A}$ and $\boldsymbol{C}$ in Lemma 1 can be replaced as follows:

$$A \to L_r, \ C \to d_k. \tag{19}$$

This is because $A$ in Lemma 1 is the arbitrary real matrix. Moreover, $d_k$ in Eq. (18) is a scalar value. Thus, the rank of $d_k$ is equal to the row number in $d_k$. As a result, $d_k$ can fulfill the condition for $C$ in Lemma 1.

Real variables $\alpha$ and $\beta$ in Lemma 1 require the matrix computation of the $A + \alpha \beta^T C$ form. Thus, using Eqs. (18) and (19), $\alpha$ and $\beta$ can be changed to $\tilde{\xi}$ and $\tilde{\psi}$ as follows:

$$\alpha \to \tilde{\xi}, \ \beta \to \tilde{\psi}. \tag{20}$$

The inequality condition $\|AC^T\|_F > 0$ in Lemma 1 can be rewritten by Eq. (19) as follows:

$$\|d_k L_r\|_F = d_k \|L_r\|_F > 0 \ \to \ d_k > 0, \|L_r\|_F > 0. \tag{21}$$

Notably, conditions $d_k > 0$ and $\|L_r\|_F > 0$ are always satisfied unless all values of $q_k$ and $L_r$ are zero. Consequently, Lemma 2 can be proven using Lemma 1. □

Based on Lemma 2, the gradients with respect to $\hat{\xi}$ and $\hat{\psi}$ are performed to derive the algebraic values of $\tilde{\xi}$ and $\tilde{\psi}$ in Eq. (18). As a result, the following formulations are required to calculate optimal solutions:

$$d_k L_r \tilde{\psi} + \tilde{\xi} \left(d_k \|\tilde{\psi}\|_2\right)^2 = \mathbf{0}, \tag{22a}$$

$$d_k L_r^T \tilde{\xi} + \left(d_k \|\tilde{\xi}\|_2\right)^2 \tilde{\psi} = \mathbf{0}, \tag{22b}$$

$$\|\tilde{\psi}\|_F \neq 0, \ \|\tilde{\xi}\|_F \neq 0. \tag{22c}$$

It should be emphasized that Eq. (22c) is automatically fulfilled to ensure inequality $\left(\|L_r + \tilde{\xi}\tilde{\psi}^T d_k\|_F\right)^2 < (\|L_r\|_F)^2$ in Lemma 2. In addition, $d_k$ in Eq. (22a) must not be zero in this derivation process, see Eq. (21). Consequently, the optimal value $\tilde{\xi}$ in Eq. (22a) can be computed without contradictions as follows:

$$\widetilde{\pmb{\xi}} = -\frac{1}{\left(d_k\|\widetilde{\pmb{\psi}}\|_2\right)^2} d_k \pmb{L}_r \widetilde{\pmb{\psi}}. \tag{23}$$

Substituting Eq. (23) into Eq. (22b), the relationship for $\widetilde{\pmb{\psi}}$ is expressed as follows:

$$\pmb{L}_r^T \pmb{L}_r \widetilde{\pmb{\psi}} = \frac{\left(\|\pmb{L}_r \widetilde{\pmb{\psi}}\|_2\right)^2}{\left(\|\widetilde{\pmb{\psi}}\|_2\right)^2} \widetilde{\pmb{\psi}}. \tag{24}$$

To comply with Eqs. (22) and (24), the optimal solution $\widetilde{\pmb{\psi}}$ in Eq (24) should be an eigenvector in the following eigenvalue problem [19-20]:

$$\pmb{L}_r^T \pmb{L}_r \widetilde{\pmb{w}} = \lambda \widetilde{\pmb{w}}, \tag{25}$$

where $\lambda$ and $\widetilde{\pmb{w}}$ denote the eigenvalue and eigenvector in Eq. (24), respectively. Note that a form of $\pmb{L}_r^T \pmb{L}_r$ in Eq. (25) is a non-zero symmetric matrix because of $\|\pmb{L}_r\|_F > 0$, see Eq. (21). Therefore, $\lambda$ and $\widetilde{\pmb{w}}$ in Eq. (25) are the real values, and eigenvalue $\lambda$ is not zero. In addition, $\pmb{L}_r^T \pmb{L}_r$ is the real scalar value because $\pmb{L}_r$ in Eq. (25) is a vector. For these reasons, the optimal solution $\widetilde{\pmb{\psi}}$ in Eq. (24), which is also eigenvector $\widetilde{\pmb{w}}$ in Eq. (25), can be determined as follows:

$$\widetilde{\pmb{\psi}} = 1 \in \mathbb{R}^1. \tag{26}$$

As a result, the optimal solution $\widetilde{\pmb{\xi}}$ in Eq. (23) is decided by substituting Eq. (26) into Eq. (23):

$$\widetilde{\pmb{\xi}} = -\frac{1}{(d_k)^2} d_k \pmb{L}_r. \tag{27}$$

The obtained solutions in Eqs. (26) and (27) can minimize the Frobenius norm in Eq. (18). Nevertheless, $\widetilde{\pmb{\xi}}$ and $\widetilde{\pmb{\psi}}$ are not auxiliary vectors $\pmb{\xi}$ and $\pmb{\psi}$ in Eq. (16). Thus, a relationship should be found between the optimal solutions in Eqs. (26) to (27) and auxiliary vectors in Eq. (16). For this reason, Lemma 3 is suggested.

Lemma 3. Assume that $\|L_r\|_F > 0$ and $d_k > 0$. Then, $\xi := L_J^T(L_J L_J^T)^{-1}\tilde{\xi}$ and $\psi := \frac{1}{d_k}q_k$ can be the auxiliary vectors that assure $\left(\|L_r + L_J \xi \psi^T q_k\|_F\right)^2 < (\|L_r\|_F)^2$.

*Proof.* The assumed $\xi$ and $\psi$ in Lemma 3 are written as follows:

$$\xi := L_J^T(L_J L_J^T)^{-1}\tilde{\xi}, \quad \psi := \frac{1}{d_k}q_k. \tag{28}$$

The optimal solution $\tilde{\psi}$ in Eq. (26) is 1. Thus, the assumed $\psi$ in Eq. (28) can be modified using $\tilde{\psi}$ as follows:

$$\psi = \frac{1}{d_k}q_k = \frac{1}{d_k}q_k \tilde{\psi}. \tag{29}$$

After substituting Eqs. (28) and (29) into Eq. (15), it is rewritten as

$$L_r + L_J \xi \psi^T q_k = L_r + L_J L_J^T(L_J L_J^T)^{-1}\tilde{\xi}\tilde{\psi}^T \left(\frac{1}{d_k}q_k\right)^T \left(d_k \frac{1}{d_k}q_k\right), \tag{30a}$$

$$\tilde{\psi} = \tilde{\psi}^T = 1 \in \mathbb{R}^1. \tag{30b}$$

Here, the magnitude of the ROM solution $q_k$ is $\|q_k\|_2 = d_k$. Thus, $\frac{1}{d_k}q_k$ is a normal vector of the ROM solution $q_k$ (i.e., $\left(\frac{1}{d_k}q_k\right)^T \left(\frac{1}{d_k}q_k\right) = 1$). Hence, Eq. (30) can be rewritten as follows:

$$L_r + L_J \xi \psi^T q_k = L_r + \tilde{\xi}\tilde{\psi}^T d_k. \tag{31}$$

Note that $\tilde{\xi}$ and $\tilde{\psi}$ in Eq. (31) can minimize $\left(\|L_r + \tilde{\xi}\tilde{\psi}^T d_k\|_F\right)^2$ under the assumption in Lemma 3, see Eqs. (21), (26), and (27). Therefore, the auxiliary vectors in Eq. (28) are the solutions in Lemma 3. □

Consequently, the additional basis $\Psi$ in the proposed method is computed by using Eqs. (14) and the relationship in Lemma 3 as follows:

$$\Psi = \Phi + \xi\psi^T = \Phi - L_J^T\left(L_J L_J^T\right)^{-1}\frac{1}{(d_k)^2}d_k L_r\left(\frac{1}{d_k}q_k\right)^T = \Phi - \frac{1}{(d_k)^2}L_J^T\left(L_J L_J^T\right)^{-1}L_r q_k^T. \quad (32)$$

Algorithm 3 describes the suggested adaptive procedure for the Newton method and Galerkin projection.

---

Algorithm 3. Proposed online adaptive basis update procedure for the Newton method and Galerkin projection

---

Input: iteration $k$, trial basis $\Phi$, ROM solver tolerance $\varepsilon_{\text{ROM}}$, FOM solver tolerance $\varepsilon_{\text{FOM}}$, the number of selection points $n_{\text{sel}}$.

Output: 1) Updated basis $[\Phi \quad \Psi]$ and 2) ROM solution $q_k$.

1. Compute the ROM solution $q_k$ using the Newton method:

    $\|\Phi^T r_k\|_2 \leq \varepsilon_{\text{ROM}}, \quad r_k = J_k \Phi q_k + F_k.$

2. Obtain the FOM error $\varepsilon$:

    $\varepsilon = \|r_k\|.$

3. If error $\varepsilon$ in Step 2 does not converge (i.e., $\varepsilon > \varepsilon_{\text{FOM}}$) do

4.     Select $n_{\text{sel}}$ row locations and then construct a Boolean matrix $B_k$ to construct the local Jacobian matrix and the local residual operator:

    $B_k = [b_1 \quad \cdots \quad b_i \quad \cdots \quad b_{n_{\text{sel}}}], \quad b_i = e_m, \quad e_m = m$ -th canonical basis, $i = 1 \ldots n_{\text{sel}},$

5.     Construct the local Jacobian matrix $L_J$ and the local residual operator $L_r$:

    $L_J = B_k^T J_k, \quad L_r = B_k^T r_k.$

6.     Calculate magnitude $d_k$ of the ROM solution $q_k$:

    $d_k = \|q_k\|_2.$

7.     Compute the additional basis $\Psi$ using the following formulation:

    $\Psi = \Phi - \frac{1}{(d_k)^2}L_J^T\left(L_J L_J^T\right)^{-1}L_r q_k^T.$

8.     Construct the updated basis $[\Phi \quad \Psi]$.

9.     Replace basis $\Phi$ using the updated basis $[\Phi \quad \Psi]$:

    $\Phi \leftarrow [\Phi \quad \Psi].$

10.    Return to Step 1.

11. end if

---

Note that the thin-QR factorization and initialization technique can be also applied to Algorithm 3 by analogy with Algorithm 2.

The additional basis $\Psi$ in Algorithm 3 requires the assumption in Eq. (21). However, Lemma 2 confirms that Eq. (21) can be automatically ensured unless all values of $L_r$ and $q_k$ are zero. In addition, the local residual operator $L_r$ is determined by freely choosing the elements of the residual operator $r_k$. Hence, the assumption in Eq. (21) can be easily secured without complicated calculations.

Let us compare F-ROM with the proposed method. The proposed method applies inverse operations in $n_{sel}$ instead of using the full model dimension $N$. Hence, the computational cost is dramatically reduced compared with F-ROM. In addition, $L_J$ and $L_r$ in Algorithm 3 are elements of $J_k$ and $r_k$ in Algorithm 2, respectively. Thus, the proposed method does not require additional preparations other than the terms in F-ROM. Consequently, the proposed method provides an efficient calculation procedure while considering the error information. The other feature is that the proposed method can add a group of vectors to the additional basis $\Psi$, unlike in the conventional approaches in which $\Psi$ is computed as a vector [17-18]. This is because the proposed algorithm acquires the additional basis via the low-rank update formulation. Accordingly, the iterative scheme in Algorithm 3 can be utilized to quickly span the high-fidelity ROM subspace. In addition, the number of the initial incremental ROM basis can be flexibly chosen. The performance of the proposed method in various nonlinear problems is numerically evaluated in the following section.

## 4. Numerical examples

In this section, we use nonlinear examples to evaluate the performance of the proposed method. Three numerical examples are covered: a Bratu's problem, a coupled viscous Burgers' equation, and a heat transfer problem of a six-pole machine [22-26]. First, the Bratu's problem is used to test the additional basis $\Psi$ of the proposed method. Next, the spanned ROM subspaces under a limited number of iterations are tested using the coupled viscous Burgers' equation. The computational costs of F-ROM and the proposed method are compared using the six-pole machine problem. Here, the Newton method and Galerkin projection are employed to obtain the numerical results. In this paper, the Boolean matrix $B_k$ in Algorithm 3 is determined by selecting $n_{sel}$ numbers in the order of magnitude of the residual operator $r_k$.

## 5.1. Bratu's problem

Let us consider the Bratu's problem, also known as the Liouville–Bratu–Gelfand problem. The Bratu's problem is widely used to demonstrate physical phenomena, such as the chemical reactor theory, solid fuel ignition model, and nanoscale applications [22]. The two-dimensional Bratu's problem is written as follows [22]:

$$\frac{\partial^2 u(x,y)}{\partial x^2} + \frac{\partial^2 u(x,y)}{\partial y^2} + \lambda e^{u(x,y)} = 0 \text{ for } \Omega = \{(x,y) \in 0 \leq x \leq 1, 0 \leq y \leq 1\} \subset \mathbb{R}^2, \quad (33a)$$

$$u(x,y) = 0 \text{ on } (x,y) \in \Gamma, \quad (33b)$$

where $x$ and $y$ denote the spatial coordinates. Moreover, the two-dimensional rectangular domain and its boundary are denoted as $\Omega$ and $\Gamma$, respectively. The positive real number $\lambda > 0$ represents the transition parameter in the Bratu's problem. The partial differential operator is expressed as $\partial$. To obtain response $u(x,y)$ in Eq. (33), the finite difference method is employed [27]. Here, $x$ and $y$ are equally divided into 50 grids.

The response $u(x,y)$ in the Bratu's problem can be divided into the following cases [22]:

$$u(x,y) = \begin{cases} \text{no solution} & (\lambda > \lambda_c) \\ \text{unique solution} & (\lambda = \lambda_c), \\ \text{two solutions (upper solution and lower solution)} & (\lambda < \lambda_c) \end{cases} \quad (34)$$

Here, $\lambda_c = 6.808124223$ is the critical value of the transition parameter. In the Bratu's problem, an initial condition $u^0(x,y)$ is determined as follows [22]:

$$u^0(x,y) = -2 \ln(2 - (1 + x - x^2)^a (1 + y - y^2)^a), \quad (35)$$

where the positive real number $a > 0$ is a parameter that determines the type of the response $u(x,y)$ in Eq. (34). The positive numbers $a$ for the lower and upper responses are chosen as 0.25 and 1.25, respectively. In addition, $a = 0.75$ is selected to take the response $u(x, y)$ at the critical transition value $\lambda_c$ [22].

To obtain the POD basis, 1000 snapshots are collected at equal intervals in $0 < \lambda \leq 2$. Next, 40 dominant vectors are selected among the POD basis vectors [9-12] to calculate the initial basis in POD-ROM and the proposed model. Here, the maximum number of modes is limited to 80, and the number of selection points $n_{sel}$ is 40. The additional basis $\Psi$ of POD-ROM is selected using the remaining POD vectors, see Eq. (11). The performance of ROM is evaluated at $\lambda = 3, 5$, and the critical transition value $\lambda = 6.808124223$.

Figs. 1, 2, and 3 represent the response fields for the full model. In addition, the differences between the full model and reduced models are also shown in Figs. 1, 2, and 3. An example procedure for the proposed method is presented in Fig. 4. These numerical results show that the proposed adaptive procedure can return high-accuracy responses, even if the initial basis cannot adequately explain the desired response fields.

*5.2. Coupled viscous Burgers' problem*

The coupled viscous Burgers' equation is widely used to describe various physics phenomena, such as wave propagation, vorticity transportation, and turbulence in hydrodynamics [23-24]. The two-dimensional coupled viscous Burgers' problem is written as follows [23-24]:

$$\frac{\partial u}{\partial t} + u\frac{\partial u}{\partial x} + v\frac{\partial u}{\partial y} = Re\left(\frac{\partial^2 u}{\partial x^2} + \frac{\partial^2 u}{\partial y^2}\right), \quad u = u(x, y, t), \tag{36a}$$

$$\frac{\partial v}{\partial t} + u\frac{\partial v}{\partial x} + v\frac{\partial v}{\partial y} = Re\left(\frac{\partial^2 v}{\partial x^2} + \frac{\partial^2 v}{\partial y^2}\right), \quad v = v(x, y, t), \tag{36b}$$

$$\Omega = \{(x, y) \in 0 \leq x \leq 1, 0 \leq y \leq 1\} \subset \mathbb{R}^2, \tag{36c}$$

where $t$ and $Re$ denote the temporal coordinate and the Reynold's number, respectively. Moreover, the velocities in $x$ and $y$ directions are denoted as $u$ and $v$, respectively. The numerical results are computed using the finite difference method [27]. Here, $x$ and $y$ are equally divided into 50 grids. In addition, the temporal coordinate $t$ is $0 \leq t \leq 1$ with a time step $\Delta t = 0.001$.

The initial and boundary conditions are determined from the following equations [23-24]:

$$u(x,y,t) = \frac{3}{4} + \frac{1}{4(1+exp[(-4x+4y-t)/(32Re)])}, \quad (37a)$$

$$v(x,y,t) = \frac{3}{4} - \frac{1}{4(1+exp[(-4x+4y-t)/(32Re)])}, \quad (37b)$$

in which $exp$ denotes the exponential function operator.

To obtain the POD basis, 1000 snapshots are collected for the Reynold number $Re = 50$ and the temporal coordinate range $0 \leq t \leq 1$. The initial basis is determined by the 5 dominant vectors among the POD basis vectors [9-12]. Moreover, the maximum mode number of POD-ROM is limited to 40, and the number of selection points $n_{sel}$ is 120. In this numerical example, we want to show that the proposed method can span the high-fidelity ROM subspace faster than F-ROM. To achieve this purpose, Algorithms 2 and 3 are used only three times when the online adaptive procedure is required. The performance of ROM is evaluated at $Re = 100$ and $Re = 1000$.

Figs. 5 and 6 represent the response fields for the full model and the difference fields in $Re = 100$. Similarly to Figs. 5 and 6, Figs. 7 and 8 show the fields in Re = 1000. An example procedure for the proposed method is illustrated in Fig. 9. According to the numerical results, the proposed method can span the accurate ROM subspace quicker than F-ROM. This is because the additional basis $\Psi$ of F-ROM is determined only as a vector and thus, the number of basis vectors may not be sufficient to explain nonlinear behavior [13]. However, unlike in the F-ROM approach, the additional basis $\Psi$ of the proposed method can be selected as a group of vectors while considering the error information. Consequently, the proposed method can be utilized to quickly supplement the high-fidelity ROM subspace.

*5.3. Heat transfer problem for the six-pole machine*

The six-pole machine is a part of rotating machines, such as the DC motor [25]. In this numerical example, we focus on the stator among the components of the six-pole machine. The temperature fields of the six-pole machine are studied to evaluate the computational procedure of the proposed method. Here, the temperature fields are numerically calculated by using the finite element method. Thus, the heat transfer problem can be formulated as follows [25,28]:

$$C\dot{T} + K\dot{T} = Q, \tag{38a}$$

$$C = \int_V \rho c H^T H \, dV, \tag{38b}$$

$$K = \int_V B^T \begin{bmatrix} k_x & 0 & 0 \\ 0 & k_y & 0 \\ 0 & 0 & k_z \end{bmatrix} B \, dV, \tag{38c}$$

in which $\rho$, $c$, $k$, and $V$ denote the density, specific heat capacity, thermal conductivity, and volume of an objective body, respectively. Vectors $T$ and $Q$ indicate the discretized temperature vector and thermal excitation, respectively. The gradient interpolation matrix and temperature interpolation matrix are expressed as $B$ and $H$, respectively. The material properties of the six-pole machine problem are selected as follows: density $\rho = 500 \frac{\text{kg}}{m^3}$, specific heat capacity $c = 200 \frac{J}{(\text{kg-K})}$, and thermal conductivity $k_x = k_y = k_z = 303.15 + 0.3T$ W/(m-K) [25-26].

Fig. 10 shows the boundary positions, observation nodes, burst random excitations, and meshing results for the six-pole machine. Specifically, temperature values in $\Gamma_T$ are maintained at 25 °C, and the excitations are applied to $\Gamma_F$. The evaluation points of ROM are nodes A and B. In this numerical example, the time range, time step, and initial temperature are set as $0 \leq t \leq 12$, $\Delta t = 0.002$, and $T_0 = 25°C$, respectively.

To obtain the POD basis, 6000 snapshots are collected using the burst excitation in $0 \leq t \leq 12$. The burst conditions for the snapshots are different from the thermal excitations in Fig. 10 except for the initial excitation and the maximum amplitude. The initial basis is determined by the 40 dominant vectors among the POD basis vectors [9-12]. Moreover, the maximum mode number is set to 200. The number of selection points $n_{sel}$ is 500. In this numerical example, the adaptive procedures are performed no more than 3 times.

Temperature results at nodes A and B are plotted in Fig. 11. In addition, Fig. 12 represents the temperature fields for the full model and difference temperature fields. The computational costs of this experiment are listed in Table 1. Furthermore, normalized relationships between the costs, iteration steps, and accuracies are illustrated in Fig. 13. Unlike the F-ROM approach, the proposed method does not require inverse operations in the full model dimension to consider the error information. Therefore, the proposed method provides an effective computational procedure, as shown in Table 1 and Fig. 13.

## 6. Conclusion

This study introduced an effective online adaptive basis technique for the reduction of nonlinear models. The proposed method can provide a highly accurate adaptive basis using the low-rank update formulation and the optimization of the residual operator. In particular, the residual operator of the proposed method is defined in local domains. Thus, the adaptive basis can be obtained with a lower computational cost compared with the conventional approaches that used full-dimensional inverse operations. In addition, this study investigates the existence conditions for the basis to avoid contradictions in the derivation process of the proposed method. As a result, the problems of costs and/or the number of basis vectors that exist in conventional approaches can be well addressed using the proposed method and its robust conditions. The performance of the proposed method is clearly confirmed in well-constructed nonlinear numerical examples. Importantly, the key idea of the study is not limited to the Galerkin projection and/or Newton method. Hence, it can be applied to other nonlinear engineering fields (e.g., the partitioning method or multi-physics problem) with various projection techniques and/or numerical solvers.

**Acknowledgments**

This research was supported by the Basic Science Research Programs of the National Research Foundation of Korea funded by the Ministry of Science, ICT, and Future Planning, South Korea (NRF-2018R1A1A1A05078730).

**References**

1. Tsompanakis, Y., Lagaros, N.D., Papadrakakis, M.: Structural design optimization considering uncertainties. Taylor & Francis, Florida, 2008

2. Nguyen, C., Zhuang, X.Y., Chamoin, L., Zhao, X.Z., Nguyen-Xuan, H., Rabczuk, T.: Three-dimensional topology optimization of auxetic metamaterial using isogeometric analysis and model order reduction. COMPUTER METHODS IN APPLIED MECHANICS AND ENGINEERING. **371** (2020)

3. Degen, D., Veroy, K., Wellmann, F.: Certified reduced basis method in geosciences Addressing the challenge of high-dimensional problems. COMPUTATIONAL GEOSCIENCES. **24**(1), 241-259 (2020)

4. Sen, S., Veroy, K., Huynh, DBP., DepariS, S., Nguyen, NC., Patera, AT.: "Natural norm" a posteriori error estimators for reduced basis approximations. JOURNAL OF COMPUTATIONAL PHYSICS. **217**(1), 37-62 (2006)

5. Giacomini, M., Veroy, K., Diez, P.: Special Issue on Credible High-Fidelity and Low-Cost Simulations in Computational Engineering. INTERNATIONAL JOURNAL FOR NUMERICAL METHODS IN ENGINEERING. **121**(23), 5151–5152 (2020)

6. Park, K.C., Park, Y.H.: Partitioned component mode synthesis via a flexibility approach. AIAA JOURNAL. **42**(6), 1236-1245 (2004)

7. Lim, J.H., Hwang, D.S., Kim, H.W., Lee, G.H., Kim, J.G.: A coupled dynamic loads analysis of satellites with an enhanced Craig-Bampton approach. AEROSPACE SCIENCE AND TECHNOLOGY. **69**, 114-122 (2017)

8. Ahn, J.G., Yang, H.I., Kim, J.G.: Multipoint constraints with Lagrange multiplier for system dynamics and its reduced-order modeling. AIAA JOURNAL. **58**(1), 385–401 (2020)

9. Amsallem, D., Farhat, C.: Interpolation method for adapting reduced-order models and application to aeroelasticity. AIAA JOURNAL. **46**(7), 1803–1813 (2008)

10. Ghavamian, F., Tiso, P., Simone, A.: POD-DEIM model order reduction for strain-softening viscoplasticity. COMPUTER METHODS IN APPLIED MECHANICS AND ENGINEERING. **317**, 458–479 (2017)

11. Kerfriden, P., Goury, O., Rabczuk, T., Bordas, S.P.A.: A partitioned model order reduction approach to rationalise computational expenses in nonlinear fracture mechanics. COMPUTER METHODS IN APPLIED MECHANICS AND ENGINEERING. **256**, 169–188 (2013)

12. Im, S., Kim, E., Cho, M.: Reduction process based on proper orthogonal decomposition for dual formulation of dynamic substructures. COMPUTATIONAL MECHANICS. **64**(5), 1237–1257 (2019)

13. Ahmed, S.E., San, O., Rasheed, A., Iliescu, T.: A long short-term memory embedding for hybrid uplifted reduced order models. PHYSICA D-NONLINEAR PHENOMENA. **409** (2020)

14. Wang, Q., Ripamonti, N., Hesthaven, JS.: Recurrent neural network closure of parametric POD-Galerkin reduced-order models based on the Mori-Zwanzig formalism. JOURNAL OF COMPUTATIONAL PHYSICS. **410** (2020)

15 Carlberg, K.: Adaptive h-refinement for reduced order models. INTERNATIONAL JOURNAL FOR NUMERICAL METHODS IN ENGINEERING. **102**(5), 1192–1210 (2015)

16. Ryckelynck, D.: A priori hyperreduction method: an adaptive approach. JOURNAL OF COMPUTATIONAL PHYSICS. **202**(1), 346-366 (2005)

17. Efendiev, Y., Gildin, E., Yang, Y.F.: Online adaptive local–global model reduction for flows in heterogeneous porous media. COMPUTATION. **4**(2) (2016)

18. Nigro, P.S.B., Anndif, M., Teixeira, Y., Pimenta, P.M., Wriggers, P.: An adaptive model order reduction with quasi-Newton method for nonlinear dynamical problems. INTERNATIONAL JOURNAL FOR NUMERICAL METHODS IN ENGINEERING. **106**(9), 740–759 (2016)

19. Peherstorfer, B., Willcox, K.: Online adaptive model reduction for nonlinear systems via low-rank updates. SIAM JOURNAL ON SCIENTIFIC COMPUTING. **37**(4), 2123–2150 (2015)

20. Golub, G. H., Van Loan, C.F.: Matrix computations, Johns Hopkins University Press, Baltimore (2012)

21. Mordecal, A.: Nonlinear programming analysis and methods, Courier Corporation, Massachusetts (2003)

22. Temimi, H., Ben-Romdhane, M., Baccouch, M., Musa, M.O.: A two-branched numerical solution of the two-dimensional Bratu's problem. APPLIED NUMERICAL MATHEMATICS. **153**, 202–216 (2020)

23. Zhang, L., Ouyang, J., Wang, XX., Zhang, XH.: Variational multiscale element-free Galerkin method for 2D Burgers' equation. JOURNAL OF COMPUTATIONAL PHYSICS. **229**(19), 7147-7161 (2010)

24. Wang, Y.P., Navon, I.M., Wang, X.Y., Cheng, Y.: 2D Burgers equation with large Reynolds number using POD/DEIM and calibration. INTERNATIONAL JOURNAL FOR NUMERICAL METHODS IN FLUIDS. **82**(12), 909–931 (2016)

25. Islam, M.S., Mir, S., Sebastian, T.: Issues in reducing the cogging torque of mass-produced permanent-magnet brushless DC motor. IEEE TRANSACTIONS ON INDUSTRY APPLICATIONS. **40**(3), 813–820 (2004)

26. Feng, S.Z., Cui, X.Y., Li, A.M.: Fast and efficient analysis of transient nonlinear heat conduction problems using combined approximations (CA) method. INTERNATIONAL JOURNAL OF HEAT AND MASS TRANSFER. **97**, 638–644 (2016)

27. LeVeque, R.J.: Finite difference methods for ordinary and partial differential equations: steady-state and time-dependent problems, Society for Industrial and Applied Mathematics, Philadelphia (2007)

28. Bathe, K.J.: Finite element procedures, Prentice Hall, New Jersey (2006)

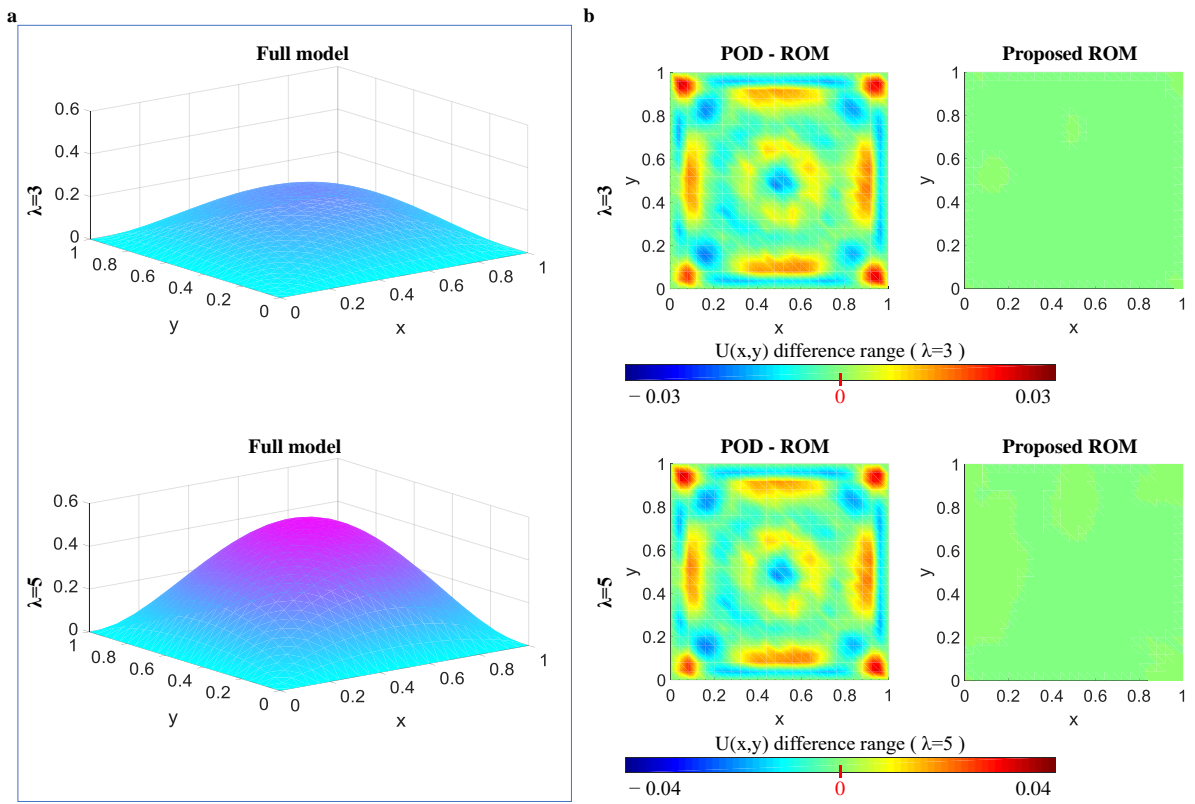

Fig. 1. Lower response fields and their differences in the Bratu's problem ($\lambda \in \{3, 5\}$): (a) lower response fields of the full model and (b) differences of the lower response fields between the full model and reduced models.

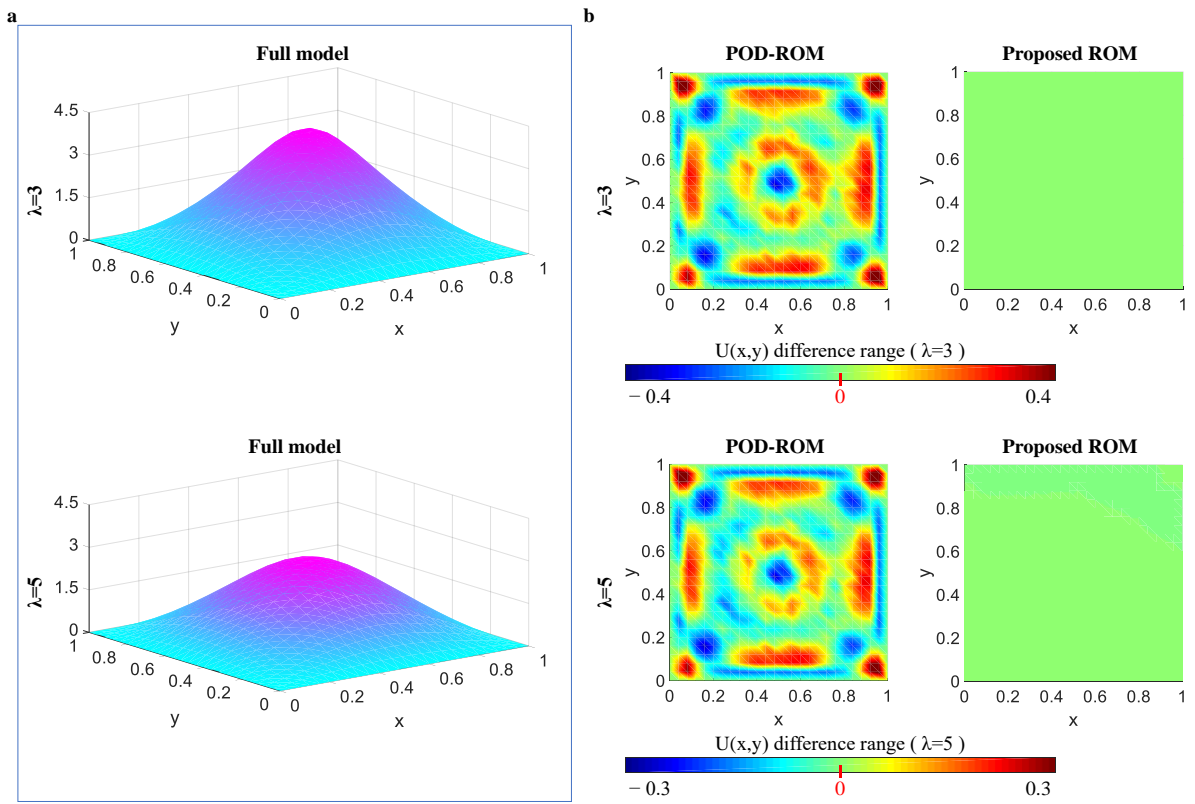

Fig. 2. Upper response fields and their differences in the Bratu's problem ($\lambda \in \{3, 5\}$): (a) upper response fields of the full model and (b) differences of the upper response fields between the full model and reduced models.

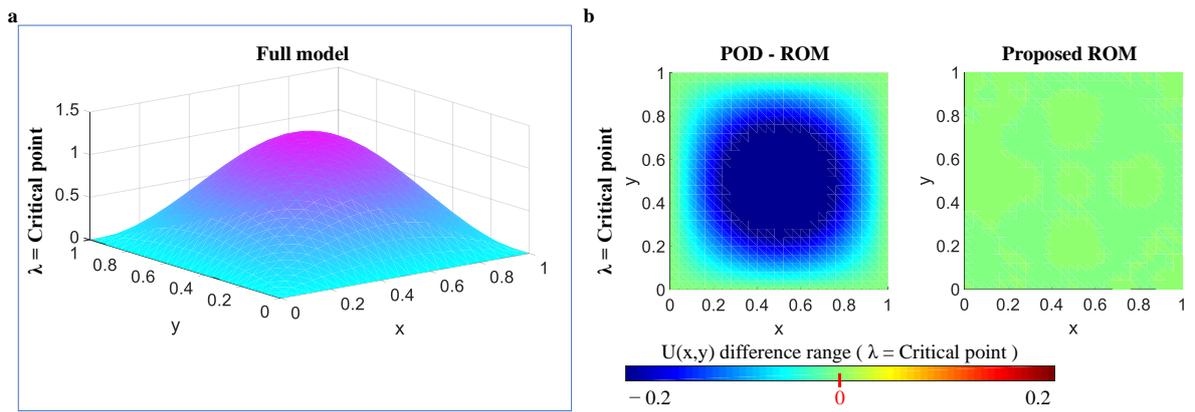

Fig. 3. Critical response field and its differences in the Bratu's problem ($\lambda = 6.8083545$): (a) critical response field of the full model and (b) differences of the critical response fields between the full model and reduced models.

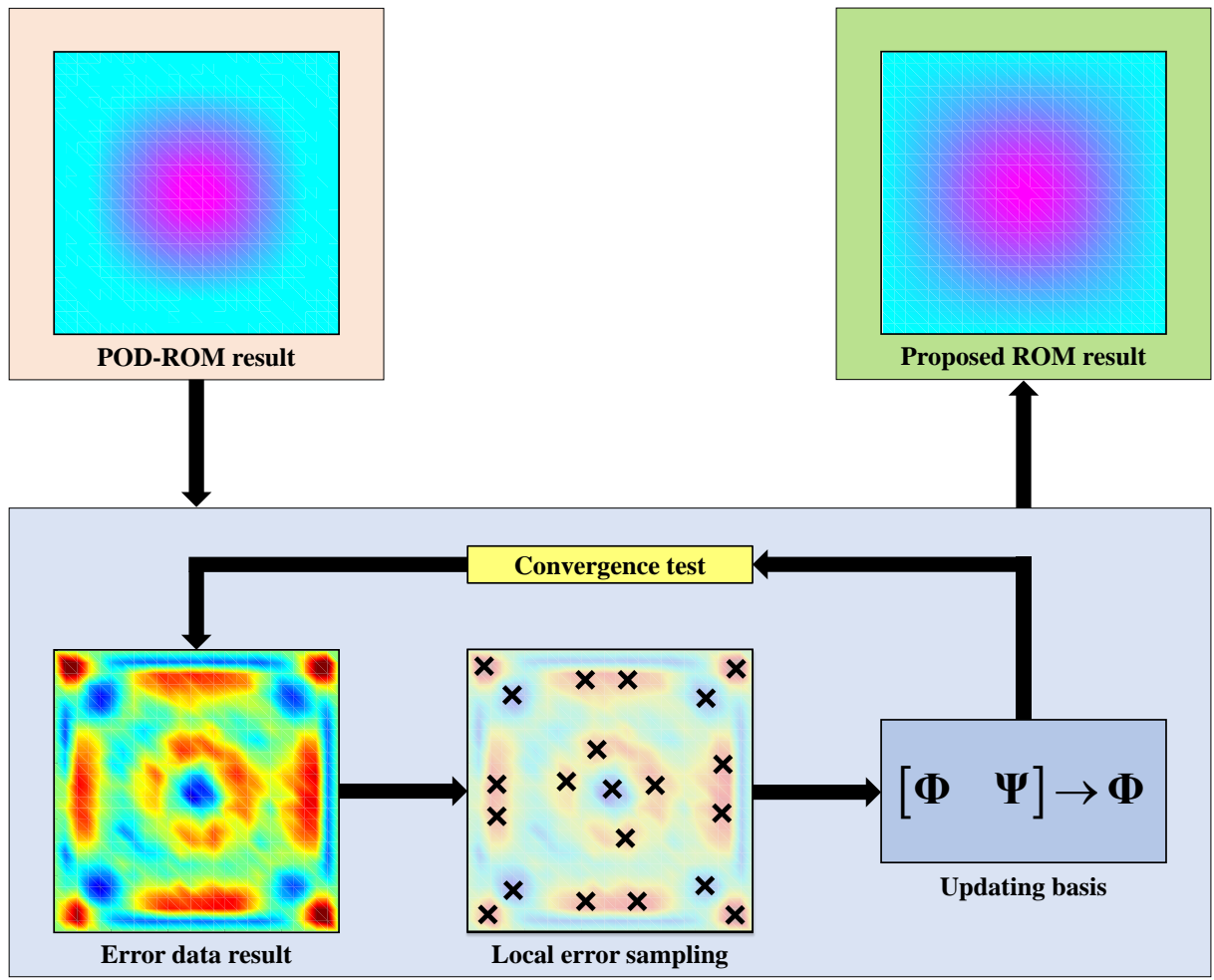

Fig. 4. Example of an application of the proposed method to the Bratu's problem.

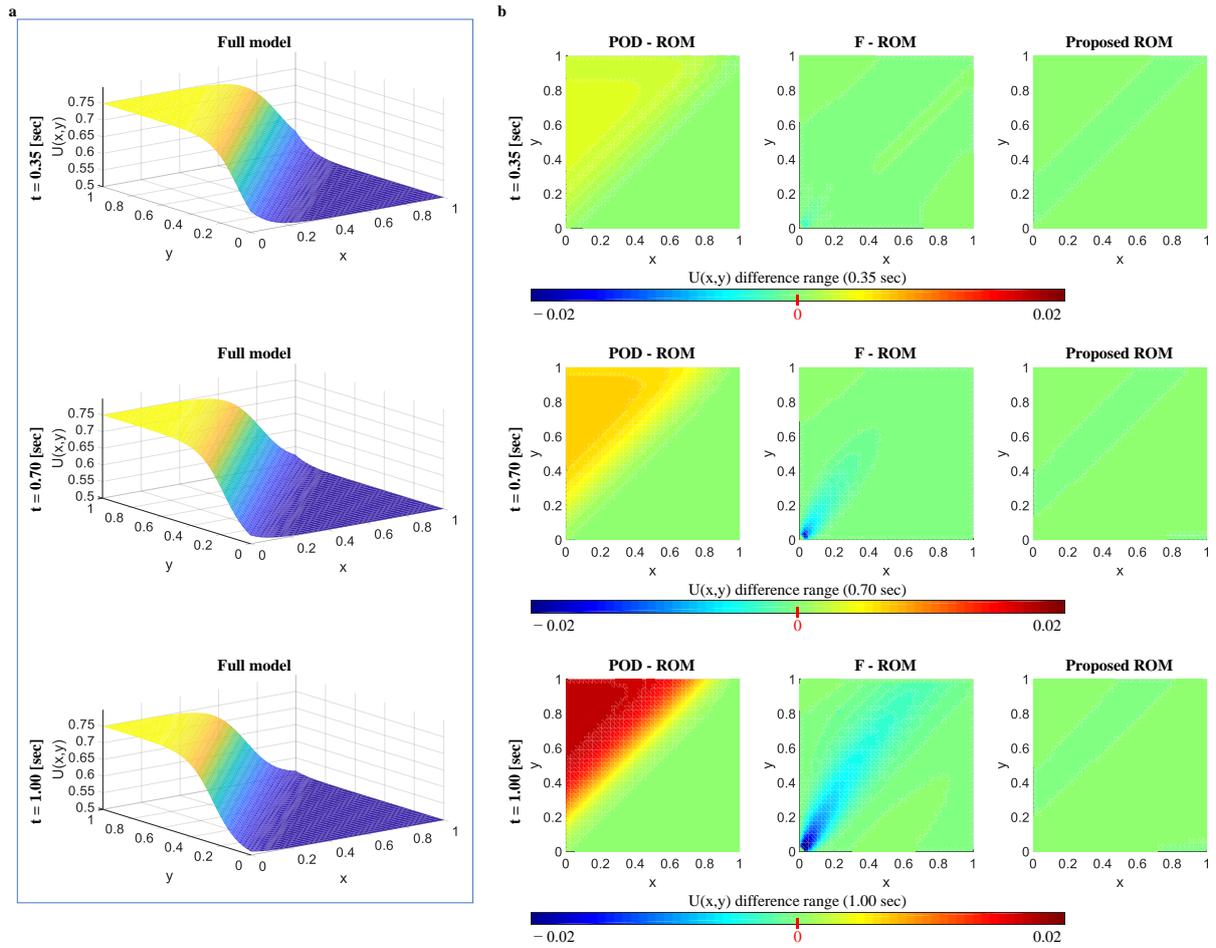

Fig. 5. $x$ direction velocity fields $u(x, y, t)$ and their differences in the coupled viscous Burgers' equation (Reynold's number = 100) : (a) velocity fields of the full model and (b) differences of the velocity fields between the full model and reduced models.

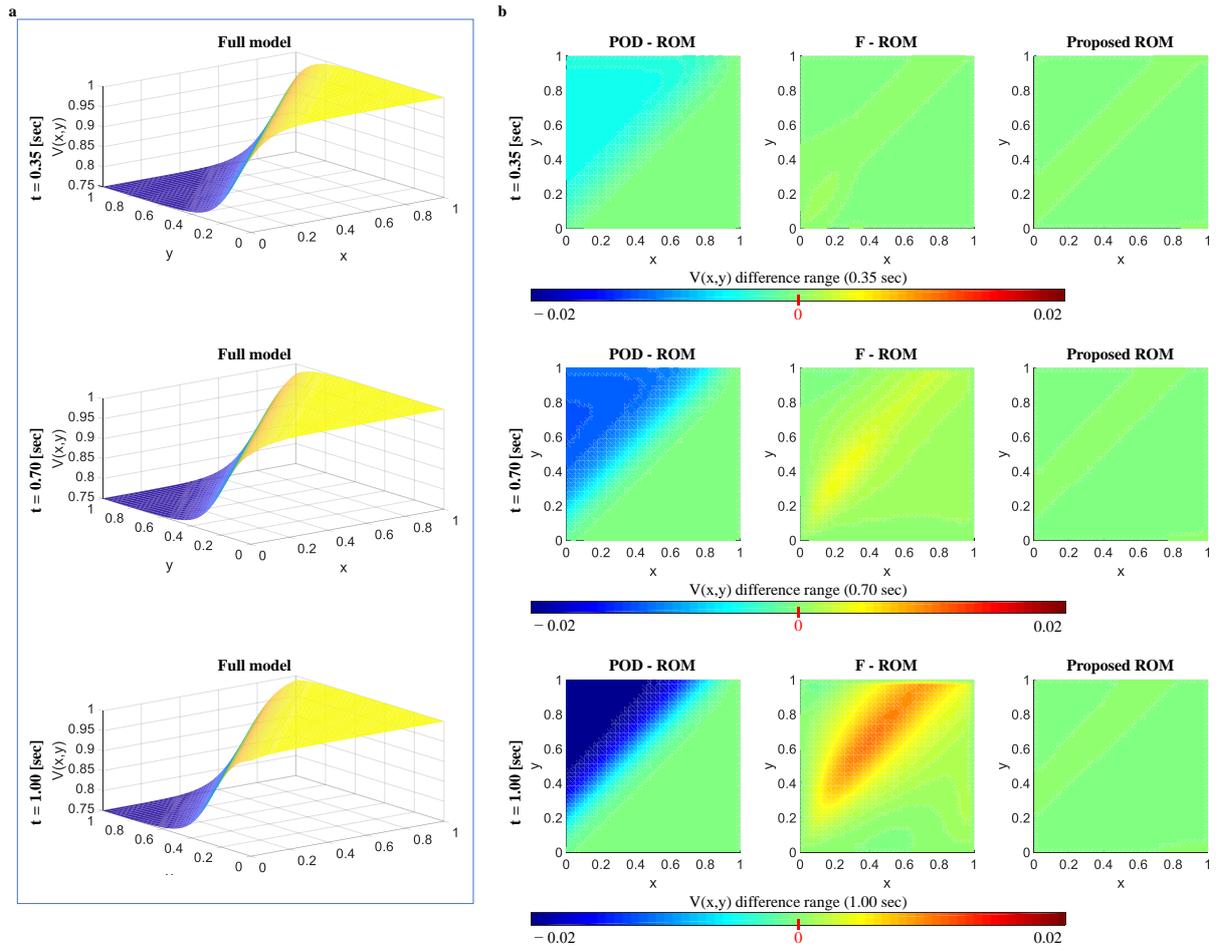

Fig. 6. $y$ direction velocity fields $v(x,y,t)$ and their differences in the coupled viscous Burgers' equation (Reynold's number = 100) : (a) velocity fields of the full model and (b) differences of the velocity fields between the full model and reduced models.

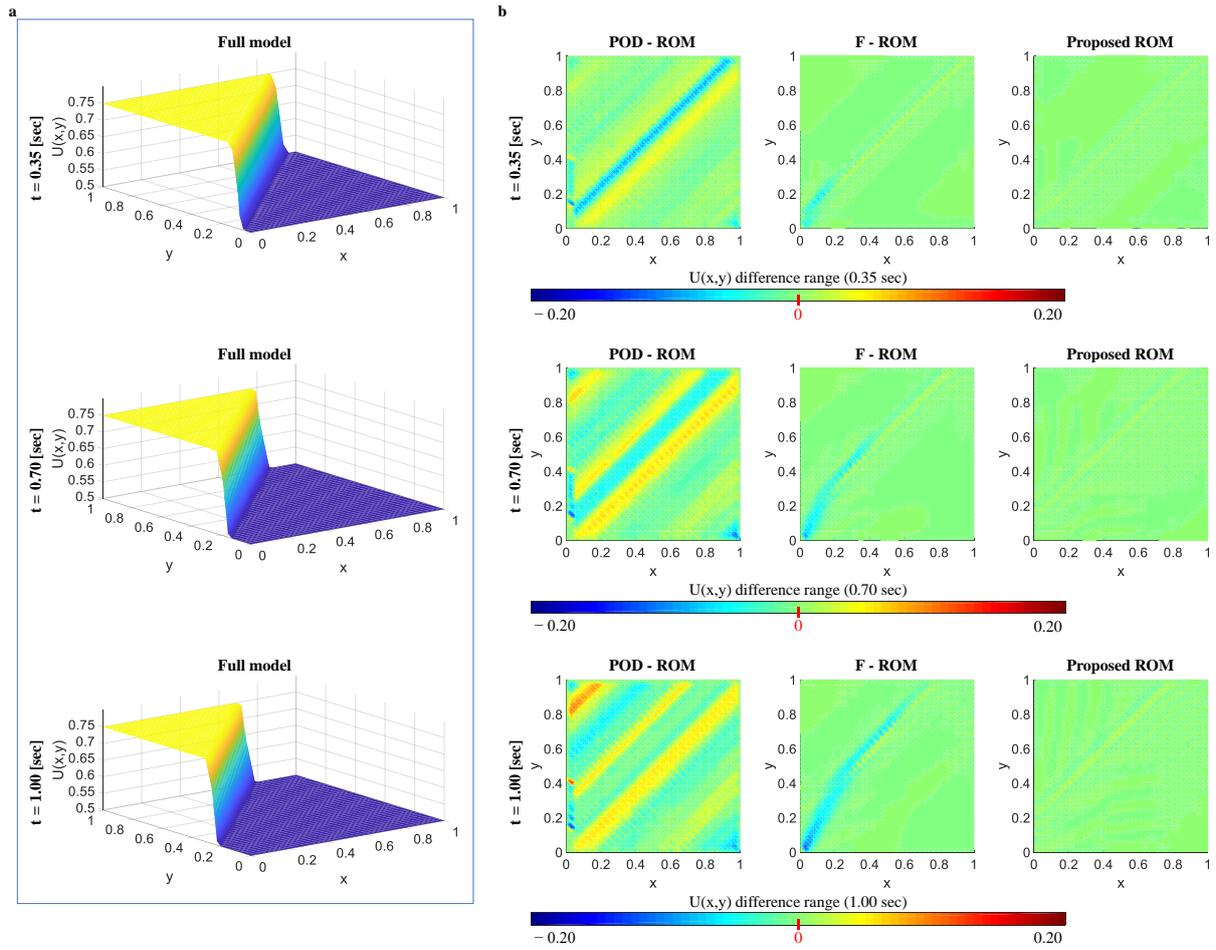

Fig. 7. $x$ direction velocity fields $u(x, y, t)$ and their differences in the coupled viscous Burgers' equation (Reynold's number = 1000) : (a) velocity fields of the full model and (b) differences of the velocity fields between the full model and reduced models.

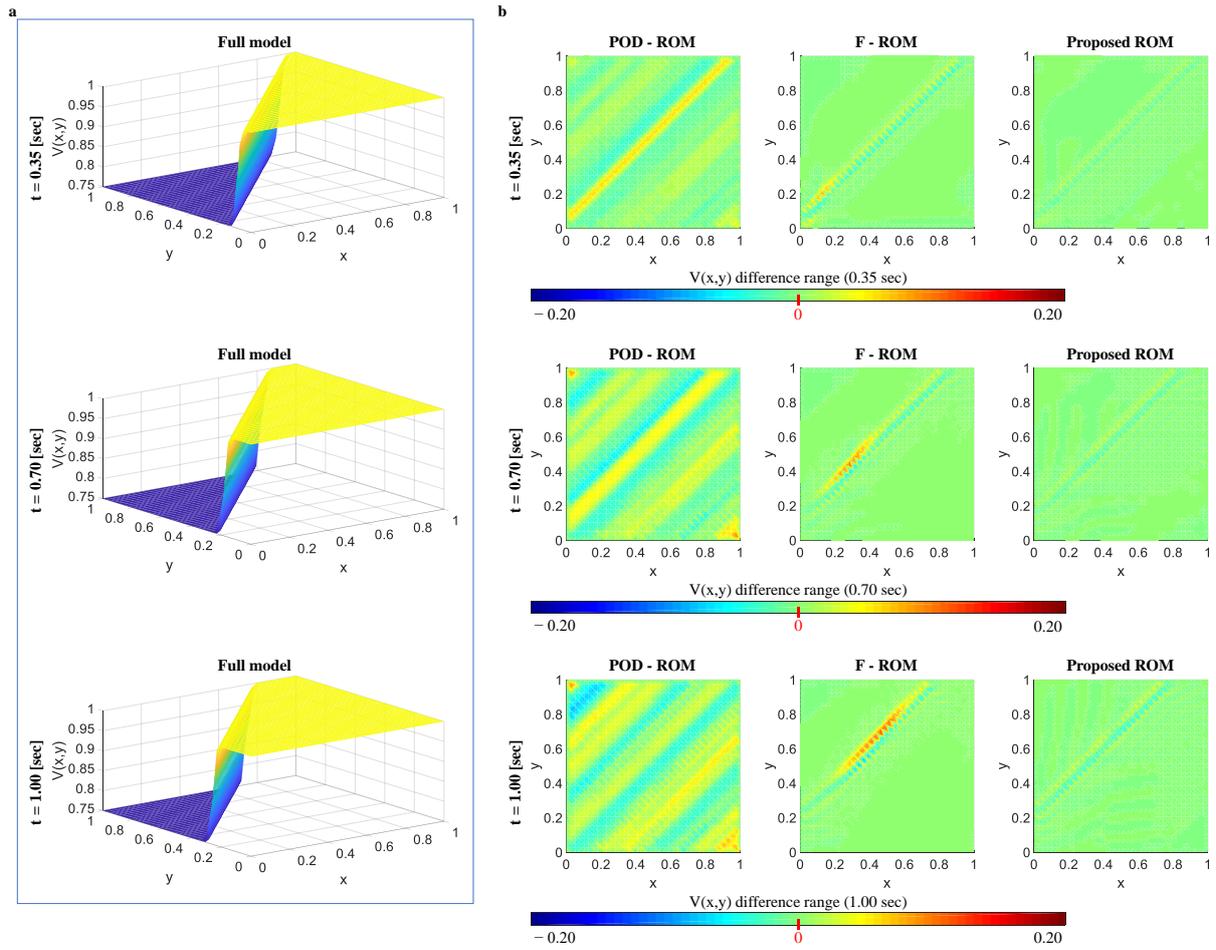

Fig. 8. $y$ direction velocity fields $v(x, y, t)$ and their differences in the coupled viscous Burgers' equation (Reynold's number = 1000) : (a) velocity fields of the full model and (b) differences of the velocity fields between the full model and reduced models.

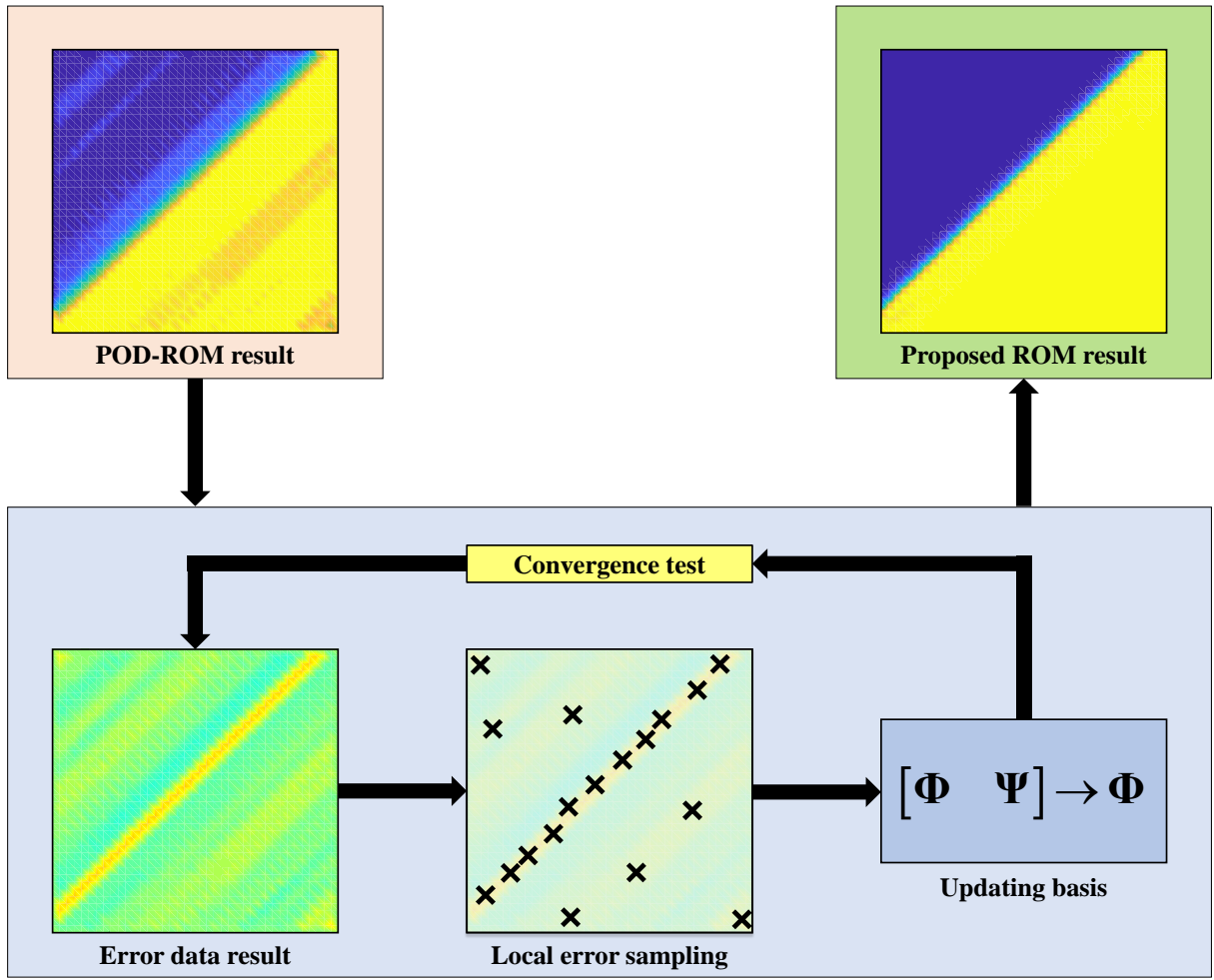

Fig. 9. Example of an application of the proposed method to the coupled viscous Burgers' problem.

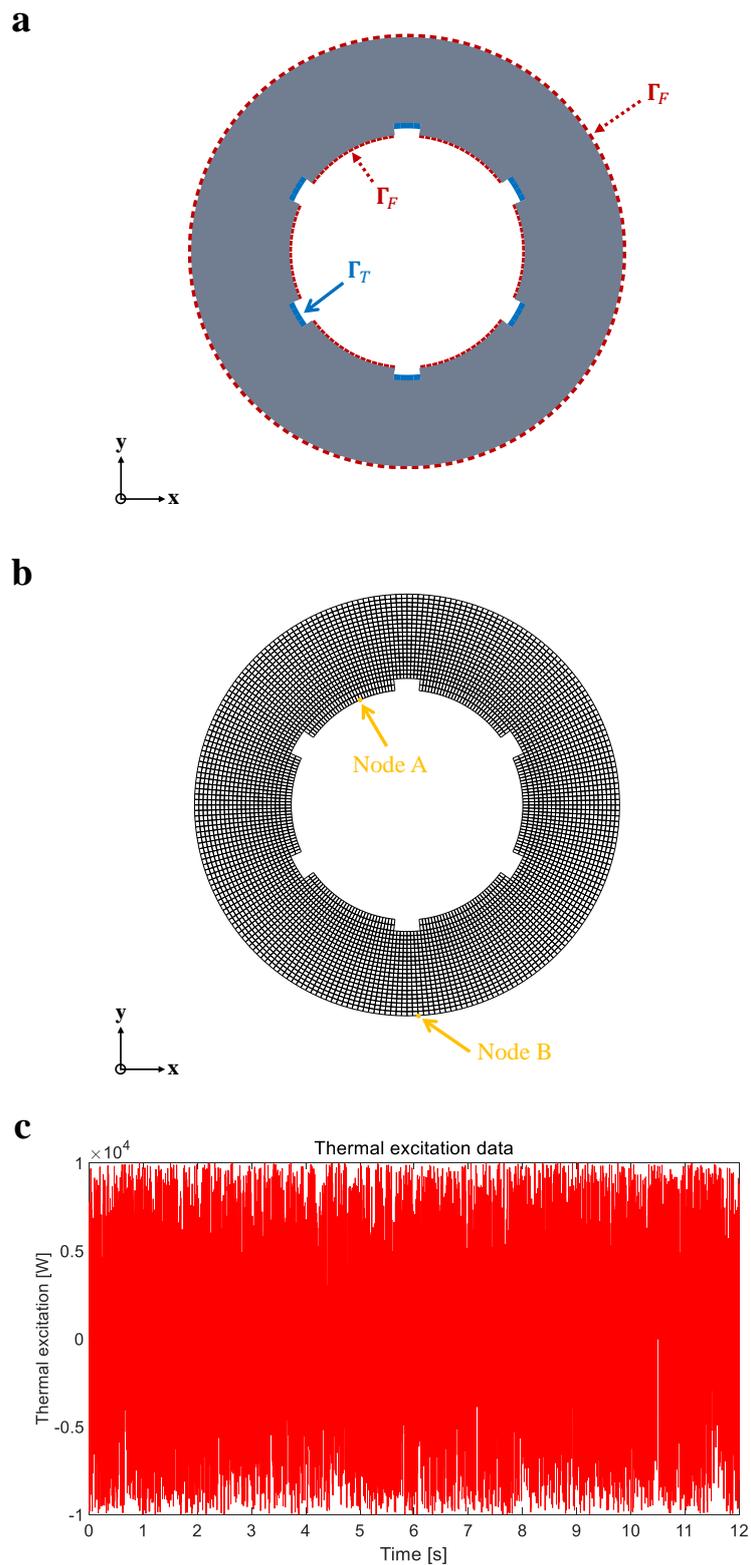

Fig. 10. Geometry, mesh result and thermal excitation conditions for the six-pole machine problem : (a) entire geometry and locations of $\Gamma_F$ and $\Gamma_T$, (b) mesh result and locations of the nodes and (c) thermal excitation data.

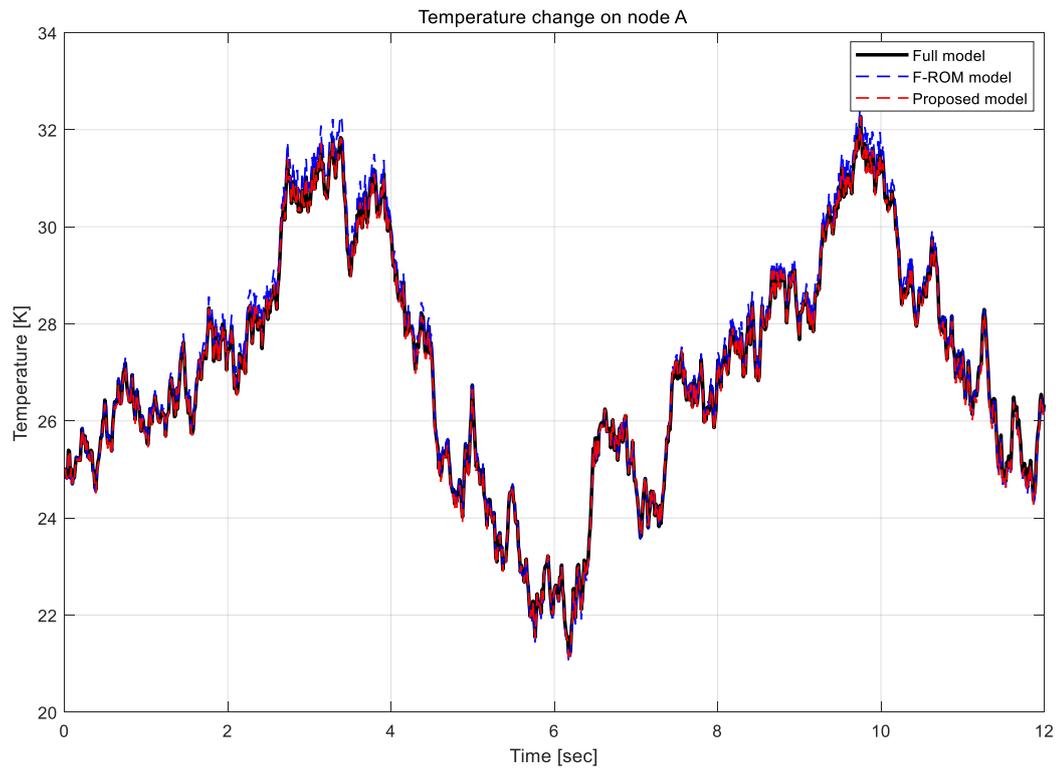

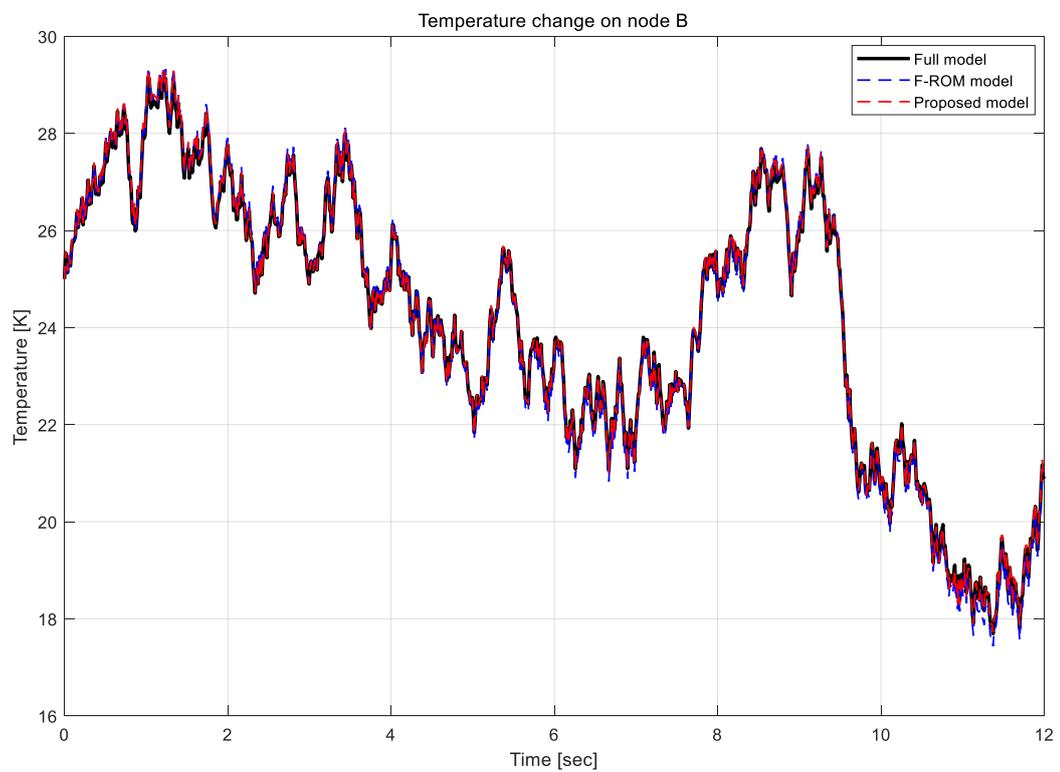

Fig. 11. Temperature results on nodes A and node B in the six pole machine problem : (a) temperature results on node A and (b) temperature results on node B.

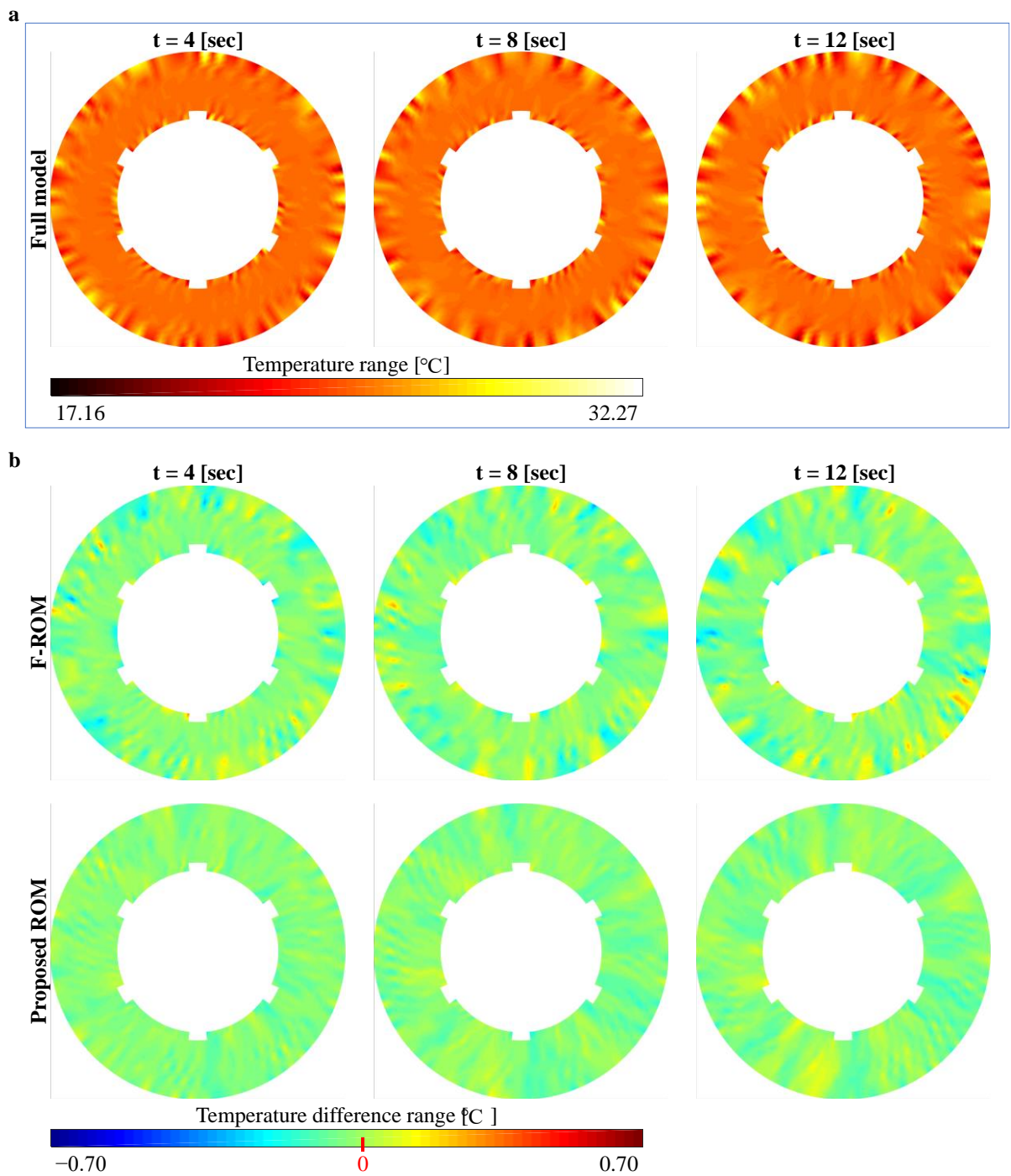

Fig. 12. Temperature fields of the six-pole machine problem: (a) temperature fields of the full model and (b) differences of the temperature fields between the full and reduced models.

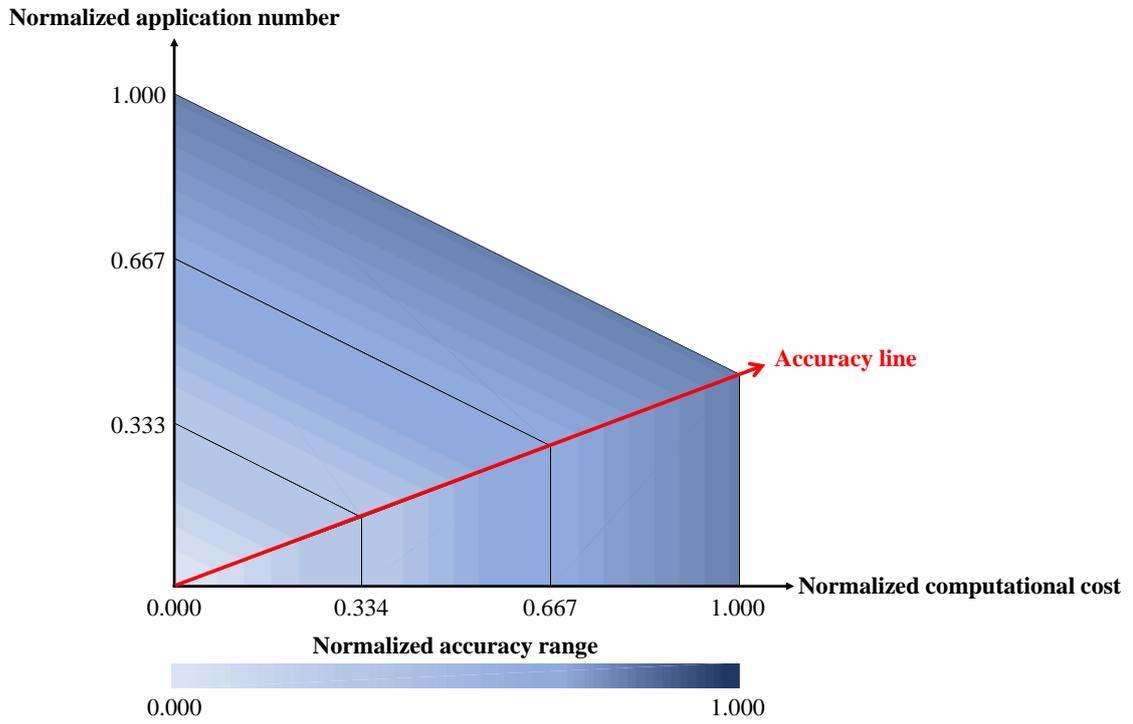

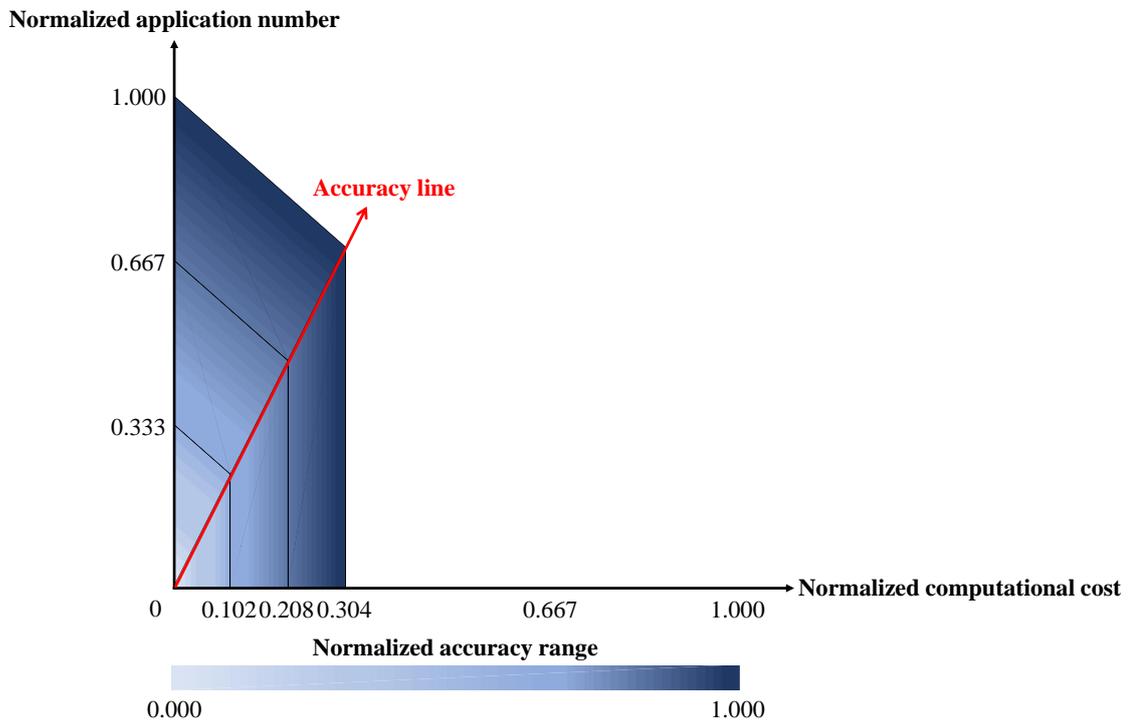

Fig. 13. Relationships between the normalized computational costs, application numbers of the adaptive procedure, and accuracies in the six-pole machine problem: (a) F-ROM and (b) proposed method.

| Selected model | Computational time [sec] |
|:---:|:---:|
| **Full model** | 78461.96 |
| **F-ROM** | 31120.74 |
| **Proposed model** | 9467.82 |

Table. 1. Computational costs in the six-pole machine problem.